\documentclass[11pt]{article}

\usepackage[square,sort,comma,numbers]{natbib}
\usepackage{fancyhdr}
\usepackage{bm,amsmath,bbm,amsfonts,amsbsy,amscd,amsxtra,amsgen,amsopn,amsthm,amssymb}
\usepackage[section]{placeins}
\usepackage{enumitem}
\usepackage{tikz,pgfplots}
\usepackage[justification=centering]{caption}
\usepackage[list=true]{subcaption}
\usepackage[colorlinks=false, pdfborder={0 0 0}]{hyperref}
\usepackage{algorithm,algorithmic}

\numberwithin{equation}{section}
\numberwithin{figure}{section}

\newcommand{\mvec}[1]{\mathrm{vec}\left(#1\right)}
\newcommand{\kron}{\otimes}
\newcommand{\trace}{\mathrm{trace}}

\newtheorem{theorem}{Theorem}[section]

\title{A low-rank approach to the solution of weak constraint variational data assimilation problems.}
\author{Melina A. Freitag\thanks{Department of Mathematical Sciences,
		University of Bath, Claverton Down, BA2 7AY, United Kingdom 
		(email: {\tt m.a.freitag@bath.ac.uk}). Corresponding author.} \and 
	Daniel L. H. Green\thanks{Department of Mathematical Sciences,
		University of Bath, Claverton Down, BA2 7AY, United Kingdom 
		(email: {\tt d.l.h.green@bath.ac.uk}).}}
\date{}

\begin{document}
	
	\maketitle
	\begin{abstract}
		\noindent Weak constraint four-dimensional variational data assimilation is an important method for incorporating data (typically observations) into a model. The linearised system arising within the minimisation process can be formulated as a saddle point problem.
		A disadvantage of this formulation is the large storage requirements involved in the linear system.
		In this paper, we present a low-rank approach which exploits the structure of the saddle point system using techniques and theory from solving large scale matrix equations.
		Numerical experiments with the linear advection-diffusion equation, and the non-linear Lorenz-95 model demonstrate the effectiveness of a low-rank Krylov subspace solver when compared to a traditional solver.
	\end{abstract}
	
	
	\section{Introduction}
	
	Data assimilation is a method for combining a numerical model with observations obtained from a physical system, in order to create a more accurate estimate for the true state of the system. One example where data assimilation is used is numerical weather prediction, however it is also applied in areas such as oceanography, glaciology and other geosciences.
	
	A property which these applications all share is the vast dimensionality of the state vectors involved. In numerical weather prediction the systems have variables of order $10^8$ \cite{Lawless2013}. In addition to the requirement that these computations to be solved quickly, the storage requirement presents an obstacle. In this paper we propose an approach for implementing the weak four-dimensional variational data assimilation method with a low-rank solution in order to achieve a reduction in storage space as well as computation time.
	The approach investigated here is based on a recent paper \cite{Stoll2015} which implemented this method in the setting of PDE-constrained optimisation. We introduce here a low-rank modification to GMRES in order to generate low-rank solutions in the setting of data assimilation.
	
	This method was motivated by recent developments in the area of solving large sparse matrix equations, see \cite{Simoncini2016, Benner2008, Saad1990, Penzl1999, Simoncini2007, Kressner2010}, notably the Lyapunov equation
	\begin{equation*}
	AX + XA^T = -BB^T
	\end{equation*}
	in which we solve for the matrix $X$, where $A$, $B$ and $X$ are large matrices of matching size. It is known that if the right hand side of these matrix equations are low-rank, there exist low-rank approximations to $X$ \cite{Grasedyck2004}. There are a number of methods which iteratively generate low-rank solutions; see e.g. \cite{Simoncini2007, Li2002, Penzl1999, Saad1990, Druskin2011}, and it is these ideas which are employed in this paper.
	
	Alternative methods \cite{Pham1998,Verlaan1997,Evensen1994} have been considered for computing low-rank solutions, based on sequential data assimilation methods such as the Kalman filter \cite{Kalman1960,Pham1998}. Furthermore there have been developments in applying traditional model reduction techniques such as Balanced Truncation \cite{Moore1981} and Principal Orthogonal Decomposition (POD) to data assimilation; e.g. \cite{Lawless2008,Cao2007}. In this paper we take a different approach, the data assimilation problem is considered in its full formulation, however the expensive solve of the linear system is done in a low-rank in time framework.
	
	In the next section we introduce a saddle point formulation of weak constraint four dimensional variational data assimilation. Section~\ref{sec:LRapproach} explains the connection between the arising linear system and the solution to matrix equations. We also introduce a low-rank approach to GMRES, and consider several preconditioning strategies.
	Numerical results are presented in Section~\ref{sec:Numerical}, with an extension to time-dependent systems considered in Section~\ref{sec:TimeDependent}.
	
	\section{Variational Data Assimilation}
	Variational data assimilation, initially proposed in \cite{Sasaki1958, Sasaki1970} is one of two families of methods for data assimilation, the other being sequential data assimilation which includes the Kalman Filter and modifications \cite{Kalman1960, Pham1998, Evensen1994}.

	We consider the discrete-time non-linear dynamical system
	\begin{equation}
	\label{eq:WeakConstraint}
	x_{k+1} = \mathcal{M}_k(x_k) + \eta_k,
	\end{equation}
	where $x_k \in \mathbb{R}^{n}$ is the state of the system at time $t_k$ and $\mathcal{M}_k: \mathbb{R}^{n} \to \mathbb{R}^{n}$ is the non-linear model operator which evolves the state from time $t_k$ to $t_{k+1}$ for $k = 0, \ldots N-1$. The model errors are denoted $\eta_k$, and are assumed to be Gaussian with zero mean and covariance matrix $Q_k \in \mathbb{R}^{n \times n}$. 
	
	Observations of this system, $y_k \in \mathbb{R}^{p_k}$ at time $t_k$ for $k = 0, \ldots N$ are given by
	\begin{equation}
	y_k = \mathcal{H}_k(x_k) + \epsilon_k,
	\end{equation}
	where $\mathcal{H}_k: \mathbb{R}^{n} \to \mathbb{R}^{p_k}$ is an observation operator, and $\epsilon_k$ is the observation error. In general, $p_k \ll n$. This observation operator $\mathcal{H}_k$ may also be non-linear, and may have explicit time dependence. The observation errors are assumed to be Gaussian, with zero mean and covariance matrix $R_k \in \mathbb{R}^{p_k \times p_k}$.
	
	We assume that at the initial time we have an a priori estimate of the state, which we refer to as the background state, and denote $x^b$. This is commonly the result of a short-range forecast, or a previous assimilation, and is typically taken to be the first guess during the assimilation process. We assume that this background state has Gaussian errors with covariance matrix $B \in \mathbb{R}^{n \times n}$.
	
	\subsection{Four dimensional variational data assimilation (4D-Var)}
	Four dimensional variational data assimilation (4D-Var) is so called for three spatial dimensions, plus time, and to differentiate it from three-dimensional variational data assimilation (3D-Var), where we do not consider multiple observation times.
	In 4D-Var, we find an initial state which minimises both the weighted least squares distance to the background state $x^b$, and the weighted least squares distance between the model trajectory of this initial state $x_k$ and the observations $y_k$ for an assimilation window $[t_0,t_N]$.
	Mathematically, we can write this as a minimisation of a cost function, e.g. $\mathrm{argmin}\; J(x)$, where
	\begin{align}
	\begin{aligned}
	\label{eq:Weak4DVarCost}
	J(x) &= \underbrace{\frac{1}{2} (x_0 - x_0^b)^T B^{-1} (x_0 - x_0^b)}_{J_b}
	+ \underbrace{\frac{1}{2} \sum_{i = 0}^{N} (y_i - \mathcal{H}_i(x_i))^T R_i^{-1} (y_i - \mathcal{H}_i(x_i))}_{J_o}
	\\&+ \underbrace{\frac{1}{2}\sum_{i=1}^{N} (x_i - \mathcal{M}_i (x_{i-1}) )^T Q_i^{-1} (x_i - \mathcal{M}_{i} (x_{i-1}) )}_{J_q}, \\
	&= \frac{1}{2} \| x_0 - x_0^b \|^2_{B^{-1}} + \frac{1}{2} \sum_{i = 0}^{N} \| y_i - \mathcal{H}_i(x_i) \|^2_{R_i^{-1}} + \frac{1}{2} \sum_{i = 1}^{N} \|x_i - \mathcal{M}_{i} (x_{i-1}) \|^2_{Q_i^{-1}},
	\end{aligned}
	\end{align}
	where $x = [x_0^T, x_1^T, \ldots, x_N^T]^T$, and $x_k$ is the model state at each timestep $t_k$ for $k = 0, \ldots, N$.
	This is known as \emph{weak constraint} 4D-Var.
	The assumption of a perfect model, gives rise to \emph{strong constraint} 4D-Var, and a simplification of the cost function, notably the removal of the $J_q$ term.
	
	The additional cost of weak constraint 4D-Var, and the difficulties in computing $Q_k$ mean that it is not widely implemented in real world systems. However, accounting for this model error (with suitable covariances) would lead to improved accuracy, and the added potential of longer assimilation windows \cite{Fisher2011a,Fisher2005}.
	
	\subsection{Incremental 4D-Var}
	To implement 4D-Var operationally, an incremental approach \cite{Courtier1994} is used, which is merely a form of Gauss-Newton iteration and generates an approximation to the solution of $x = \mathrm{argmin}\; J(x)$.
	We approximate the 4D-Var cost function by a quadratic function of an increment ${\delta x^{(\ell)}} = \left[(\delta x_0^{(\ell)})^T, (\delta x_1^{(\ell)})^T, \ldots, (\delta x_N^{(\ell)})^T\right]^T$ defined as
	\begin{equation}
	\delta x^{(\ell)} = x^{(\ell+1)} - x^{(\ell)},
	\end{equation}
	where ${x^{(\ell)}} = \left[(x_0^{(\ell)})^T, (x_1^{(\ell)})^T, \ldots, (x_N^{(\ell)})^T\right]^T$ denotes the $\ell$-th iterate of the Gauss-Newton algorithm. Updating this estimate is implemented in an \emph{outer loop}, whilst generating $\delta x^{(\ell)}$ is referred to as the \emph{inner loop}.
	This increment $\delta x^{(\ell)}$ is a solution to the minimisation of the linearised cost function 
	\begin{align}
	\begin{aligned}
	\label{eq:WeakInc}
	\tilde{J}(\delta x^{(\ell)}) &= 
	\frac{1}{2} (\delta x_0^{(\ell)} - b_0^{(\ell)})^T B^{-1} (\delta x_0^{(\ell)} - b_0^{(\ell)})\\
	&+ \frac{1}{2} \sum_{i = 0}^{N} (d_i^{(\ell)} - H_i \delta x_i^{(\ell)})^T R_i^{-1} (d_i^{(\ell)} - H_i \delta x_i^{(\ell)})\\
	&+ \frac{1}{2}\sum_{i=1}^{N} (\delta x_i^{(\ell)} - M_i \delta x_{i-1}^{(\ell)} - c_k^{(\ell)})^T Q_i^{-1} (\delta x_i^{(\ell)} - M_i \delta x_{i-1}^{(\ell)} - c_k^{(\ell)}).
	\end{aligned}
	\end{align}
	Here $M_k \in \mathbb{R}^{n \times n}$ and $H_k \in \mathbb{R}^{n \times p_k}$, are linearisations of $\mathcal{M}_k$ and $\mathcal{H}_k$ about the current state trajectory $x^{(\ell)}$. For convenience and conciseness, we introduce
	\begin{align}
	b_0^{(\ell)} &= x_0^b - x_0^{(\ell)}, \\
	d_k^{(\ell)} &= y_k - \mathcal{H}_k(x_k^{(\ell)}), \\
	c_k^{(\ell)} &= \mathcal{M}_{k} (x_{k-1}^{(\ell)}) - x_k^{(\ell)}.
	\end{align}
	
	We define the following vectors in order to rewrite the cost function in a more compact form.
	\begin{equation*}
	\delta x = \begin{bmatrix} \delta x_0 \\ \delta x_1 \\ \vdots \\ \delta x_N \end{bmatrix}, \quad
	\delta p = \begin{bmatrix} \delta x_0 \\ \delta q_1 \\ \vdots \\ \delta q_N \end{bmatrix},
	\end{equation*}
	where we have dropped the superscript for the outer loop iteration. These two vectors are related by $\delta q_k = \delta x_k - M_k \delta x_{k-1}$, or in matrix form
	\begin{equation}
	\delta p = L \delta x,
	\end{equation}
	where \begin{equation}
	\label{eq:L}
	L = \begin{bmatrix}
	I 	 &  	  & 	   & \\
	-M_1 & I  	  & 	   & \\
	& \ddots & \ddots & \\
	&		  & -M_N   & I
	\end{bmatrix} \in \mathbb{R}^{(N+1)n \times (N+1)n}.
	\end{equation}
	Furthermore, we introduce the following matrices:
	\begin{equation*}
	D = \begin{bmatrix}
	B &&& \\ &Q_1&& \\ &&\ddots& \\ &&&Q_N
	\end{bmatrix} \in \mathbb{R}^{(N+1)n \times (N+1)n},\quad
	\mathcal{R} = \begin{bmatrix}
	R_0 &&& \\ &R_1&& \\ &&\ddots& \\ &&&R_N
	\end{bmatrix} \in \mathbb{R}^{\sum\limits_{k=0}^{N}p_k \times \sum\limits_{k=0}^{N}p_k},
	\end{equation*}
	\begin{equation*}
	\mathcal{H} = \begin{bmatrix}
	H_0 &&& \\ &H_1&& \\ &&\ddots& \\ &&&H_N
	\end{bmatrix} \in \mathbb{R}^{(N+1)n \times \sum\limits_{k=0}^{N}p_k}, 
	b = \begin{bmatrix}
	b_0 \\ c_1 \\ \vdots \\ c_N
	\end{bmatrix} \in \mathbb{R}^{(N+1)n},\quad
	d = \begin{bmatrix}
	d_0 \\ d_1 \\ \vdots \\ d_N
	\end{bmatrix} \in \mathbb{R}^{\sum\limits_{k=0}^{N}p_k}.
	\end{equation*}
	
	This allows us to write \eqref{eq:WeakInc}, with the superscripts dropped, as a function of $\delta x$:
	\begin{equation}
	\tilde{J}(\delta x) = (L \delta x - b)^T D^{-1} (L \delta x - b) + (\mathcal{H} \delta x - d)^T \mathcal{R}^{-1} (\mathcal{H} \delta x - d).
	\end{equation}
	
	Minimising the cost function is equivalent to solving the linear system for the gradient. Indeed, taking the gradient of this cost function with respect to $\delta x$, we have
	\begin{equation}
	\nabla \tilde{J}(\delta x) = L^T D^{-1} (L \delta x - b) + \mathcal{H}^T \mathcal{R}^{-1} (\mathcal{H} \delta x - d).
	\end{equation}
	
	Defining $\lambda = D^{-1}(b - L \delta x)$ and $\mu = \mathcal{R}^{-1} (d - \mathcal{H} \delta x)$, allows us to write the gradient at the minimum as 
	\begin{align}
	\nabla \tilde{J} = L^T\lambda + H^T \mu &= 0. \label{eq:SaddleJ}
	\intertext{Additionally, we have}
	D\lambda + L \delta x &= b, \label{eq:Saddlelambda}\\
	R\mu + H \delta x &=d, \label{eq:Saddlemu}
	\end{align}
	and \eqref{eq:SaddleJ}, \eqref{eq:Saddlelambda} and \eqref{eq:Saddlemu} can be combined into a single linear system:
	\begin{equation}
	\label{eq:WeakSaddle}
	\begin{bmatrix}
	D & 0 & L \\
	0 & \mathcal{R} & \mathcal{H} \\
	L^T & \mathcal{H}^T & 0
	\end{bmatrix} \begin{bmatrix}
	\lambda \\ \mu \\ \delta x
	\end{bmatrix} = \begin{bmatrix}
	b \\ d \\ 0
	\end{bmatrix},
	\end{equation}
	which is solved for $\delta x$.
	
	This equation is known as the saddle-point formulation for weak constraint 4D-Var, and allows us to exploit the saddle point structure for linear solves and preconditioning \cite{Benzi2005, Stoll2015, Bergamaschi2007}.
	
	The saddle point matrix in \eqref{eq:WeakSaddle}, is a square symmetric indefinite matrix of size $\left(2n(N+1) + \sum_{k=0}^{N}p_k \right)$. In order to successfully solve this system we must use an iterative solver such as MINRES or GMRES as it is unfeasible with these large problem sizes to use a direct method. Additionally we require a good choice of preconditioner for a saddle point system \cite{Benzi2005,Benzi2008,Bergamaschi2007, Bergamaschi2009,Bergamaschi2011,Fisher2011a}, which in a data assimilation setting, has a $(1,2)$ block which is more computationally expensive than the $(1,1)$ block. The inexact constraint preconditioner \cite{Bergamaschi2007} has been found to be an effective choice of preconditioner for the data assimilation problem \cite{Fisher2011a}, but application of this results in a nonsymmetric system necessitating the use of GMRES. We consider different preconditioning approaches in Section~\ref{sec:Preconditioning}.
	Furthermore, to overcome the storage requirements of the matrix in \eqref{eq:WeakSaddle}, we wish to avoid forming it (and indeed as many of the submatrices as possible), which motivates the method described in the following section.
	
	\section{Low-rank approach}
	\label{sec:LRapproach}
	\subsection{Kronecker formulation}
	\label{sec:Kronecker}
	As noted above, the matrix formed in the saddle point formulation is very large, as indeed are the vectors $\lambda, \mu, \delta x$. We wish to adapt the ideas developed in \cite{Stoll2015} in order to solve \eqref{eq:WeakSaddle}. This approach is dependent on the Kronecker product and the $\mvec{\cdot}$ operator; which are defined to be
	\begin{equation*}
	\mathcal{A} \kron \mathcal{B} = \begin{bmatrix} a_{11} \mathcal{B} & \cdots & a_{1n} \mathcal{B} \\ \vdots & \ddots & \vdots \\ a_{m1} \mathcal{B} & \cdots & a_{mn} \mathcal{B} \end{bmatrix} \quad 
	\mvec{\mathcal{C}} = \begin{bmatrix} c_{11} \\ \vdots \\ c_{1n} \\ \vdots \\ c_{mn}\end{bmatrix}.
	\end{equation*}
	We also make use of the relationship between the two: 
	\begin{equation}
	\label{eq:SylKronRelation}
	(\mathcal{B}^T \kron \mathcal{A})\mvec{\mathcal{C}} = \mvec{\mathcal{ACB}}.
	\end{equation}
	
	Employing this definition, we may rewrite \eqref{eq:WeakSaddle} as
	\begin{equation}
	\label{eq:KroneckerWeakSaddle}
	\begin{bmatrix}
	E_1 \kron B + E_2 \kron Q & 0 & I_{N+1} \kron I_n + C \kron M \\
	0 & I_{N+1} \kron R & I_{N+1} \kron H \\
	I_{N+1} \kron I_n + C^T \kron M^T & I_{N+1} \kron H^T & 0
	\end{bmatrix} \begin{bmatrix}
	\lambda \\ \mu \\ \delta x
	\end{bmatrix} = \begin{bmatrix}
	b \\ d \\ 0
	\end{bmatrix},
	\end{equation}
	where we make the additional assumptions that $Q_i = Q$, $R_i = R$, $H_i = H$, $M_i = M$ and the number of observations $p_i = p$ for each $i$. The extended case relaxing this assumption is considered in Section~\ref{sec:TimeDependent}.
	Here 
	\begin{equation*}
	C = \begin{bmatrix}
	0&&&\\
	-1&0&&\\
	&\ddots&\ddots& \\
	&&-1&0
	\end{bmatrix},\quad E_1 = \begin{bmatrix}
	1&&&\\
	&0&&\\
	&&\ddots& \\
	&&&0
	\end{bmatrix}, \;\; \text{and} \; E_2 = \begin{bmatrix}
	0&&&\\
	&1&&\\
	&&\ddots& \\
	&&&1
	\end{bmatrix}.
	\end{equation*} 
	The matrices $C, E_1, E_2, I_{N+1} \in \mathbb{R}^{N+1 \times N+1}$, whilst $B, Q, M, I_n \in \mathbb{R}^{n \times n}, H \in \mathbb{R}^{p \times n}$, and $R \in \mathbb{R}^{p \times p}$.
	
	Using \eqref{eq:SylKronRelation}, we may rewrite \eqref{eq:KroneckerWeakSaddle} as the simultaneous matrix equations:
	\begin{align}
	\label{eq:MatrixWeakSaddle}
	\begin{aligned}
	B \Lambda E_1 + Q \Lambda E_2 + X + M X C^T &= \mathbbm{b}, \\
	R U + H X &= \mathbbm{d}, \\
	\Lambda + M^T \Lambda C + H^T U &= 0.,
	\end{aligned}
	\end{align}
	where we suppose $\lambda, \delta x, b, \mu$ and $d$ are vectorised forms of the matrices $\Lambda,X,\mathbbm{b}  \in \mathbb{R}^{n \times N+1}$ and $U, \mathbbm{d}\in \mathbb{R}^{p \times N+1}$ respectively.
	These are generalised Sylvester equations, which we solve for $\Lambda, U$ and $X$, though for implementing incremental data assimilation, we require only $\delta x$ and hence the solution $X$.
	
	For standard Sylvester equations of the form $\mathcal{AX} + \mathcal{XB} = \mathcal{C}$, it is known that if the right hand side $\mathcal{C}$ is low-rank, then there exist low-rank approximate solutions \cite{Grasedyck2004}.
	Indeed, recent algorithms for solving these Sylvester equations have focused on constructing low-rank approximate solutions. These algorithms include Krylov subspace methods (see \cite{Simoncini2016}) and ADI based methods (see \cite{Benner2009,Benner2014,Flagg2013}).
	It is this knowledge which motivates the following approach.
	
	\subsection{Existence of a low-rank solution}
	\label{sec:Existence}
	
	We wish to show that we can find a low-rank approximate solution to \eqref{eq:KroneckerWeakSaddle}. Further to the assumption that the model and observations are not time-dependent, let us additionally assume that the model is linear and perfect. Thus $c_k = M_k(x_{k-1}) - x_k = 0$ for all $k$, giving
	\begin{equation}
	\label{eq:Low-Rankb}
	b = \begin{bmatrix}
	b_0 \\ c_1 \\ \vdots \\ c_N
	\end{bmatrix} 
	= \begin{bmatrix}
	b_0 \\ 0 \\ \vdots \\ 0
	\end{bmatrix}, \qquad \text{and hence } \mathbbm{b} = \begin{bmatrix}
	b_0 & 0 & \cdots & 0
	\end{bmatrix} \in \mathbb{R}^{n \times N+1}.
	\end{equation}
	
	Assuming $R$ is non-singular, solving the second block-row of \eqref{eq:KroneckerWeakSaddle} for $\mu$ yields, 
	\begin{equation}
	\mu = (I_{N+1} \kron R^{-1}) d - (I_{N+1} \kron R^{-1}H)\delta x,
	\end{equation}
	which when substituted into the third block-row of \eqref{eq:KroneckerWeakSaddle} gives
	\begin{equation}
	(I_{N+1} \kron I_n + C^T \kron M^T) \lambda - (I_{N+1} \kron H^TR^{-1}H)\delta x = - (I_{N+1} \kron H^TR^{-1}) d.
	\end{equation}
	Reformulating this as a matrix equation as before, we are left with the simultaneous (block-row) equations
	\begin{align}
	X + MXC^T + B \Lambda E_1 + Q \Lambda E_2 &= \mathbbm{b} \label{eq:sim1}\\
	\Lambda + M^T \Lambda C - H^TR^{-1}HX &= -H^TR^{-1}\mathbbm{d}. \label{eq:sim2}
	\end{align}
	Assuming $M^{-1}$ exists, we multiply \eqref{eq:sim2} by $M^{-T}$ to obtain
	\begin{equation}
	M^{-T} \Lambda + \Lambda C = M^{-T}H^TR^{-1}HX -M^{-T}H^TR^{-1}\mathbbm{d}. \label{eq:SylLambda}
	\end{equation}
	Typically in real world applications, we only observe a small proportion of the state space. As such, the matrix $\mathbbm{d}$ containing these observations is low-rank, as is the observation operator $H$. Hence the right hand side of \eqref{eq:SylLambda} is low-rank.
	
	Applying the existence of low-rank solutions for Sylvester equations shown in \cite{Grasedyck2004} to \eqref{eq:SylLambda}, we have that $\Lambda$, or indeed an approximate solution $\tilde{\Lambda}$, is low-rank.
	
	Finally, multiplying \eqref{eq:sim1} by $M^{-1}$, and substituting in $\tilde{\Lambda}$ gives another Sylvester equation of the form
	\begin{equation}
	\label{eq:SylX}
	M^{-1}X + XC^T = M^{-1}\left(\mathbbm{b} - B \tilde{\Lambda} E_1 - Q \tilde{\Lambda} E_2\right).
	\end{equation}
	
	From the assumption that the model is perfect, we see from \eqref{eq:Low-Rankb} that $\mathbbm{b}$ is indeed low-rank, being rank 1, and hence from above, so is $\tilde{\Lambda}$.
	Thus the right hand side of this Sylvester equation \eqref{eq:SylX} is also low-rank. Applying once more the result from \cite{Grasedyck2004}, we obtain the desired property that $X$ is low-rank, or indeed there is an approximate solution $\tilde{X}$ to $X$ which is low-rank.
	
	We formulate this result as the following Theorem.
	
	\begin{theorem}
		Consider the solution to the saddle point formulation of the linearised weak constraint 4D-Var problem \eqref{eq:MatrixWeakSaddle}. Let the model and observations be time-independent, with $M = M_k, R = R_k, H = H_k, Q = Q_k$ for all $k$. Furthermore, we assume there is no model error, and that the model operator $M$, and the covariance matrix $R$ are invertible. If the number of observations $p \ll n$, then there exists a low-rank approximation $X_r = WV^T$ to $X$, where $\delta x = \mvec{X}$.
	\end{theorem}

	It is necessary to note that it would be unfeasible to compute low-rank solutions to \eqref{eq:WeakSaddle} in such a way. Indeed in \eqref{eq:SylLambda} the right hand side still contains $X$, however the observation operator allows us to know the right hand side is low-rank.
	
	Furthermore we had to make a number of assumptions to obtain this result. Whilst the assumption that $m \ll n$ is realistic, the constant operators and covariance matrices are restrictive. However, as we will see in Section~\ref{sec:TimeDependent}, relaxing some of these assumptions still results in low-rank solutions observed numerically.

	\subsection{Low-Rank GMRES (LR-GMRES)}
	\label{sec:LR-GMRES}
	In order to implement the above, we suppose as in \cite{Stoll2015,Benner2013}, that the matrices $\Lambda, U, X$ in \eqref{eq:MatrixWeakSaddle} have low-rank representations, with 
	\begin{align}
	\Lambda &= W_{\Lambda}V_{\Lambda}^T, \qquad W_{\Lambda} \in \mathbb{R}^{n \times k_{\Lambda}}, V_{\Lambda} \in \mathbb{R}^{N+1 \times k_{\Lambda}}, \\
	U &= W_UV_U^T, \qquad W_U \in \mathbb{R}^{p \times k_U}, V_U \in \mathbb{R}^{N+1 \times k_U}, \\
	X &= W_XV_X^T, \qquad W_X \in \mathbb{R}^{n \times k_X}, V_X \in \mathbb{R}^{N+1 \times k_X},
	\end{align}
	where $k_{\Lambda}, k_U, k_X \ll n$ and $k_{\Lambda}, k_U, k_X \ll N$.
	
	This allows us to rewrite \eqref{eq:MatrixWeakSaddle} as follows:
	\begin{align}
	\begin{aligned}
	\label{eq:LRMatrixWeakSaddle}
	\begin{bmatrix}
	BW_{\Lambda} & QW_{\Lambda} & W_{X} & MW_{X}
	\end{bmatrix}
	\begin{bmatrix}
	V_{\Lambda}^TE_1 \\ V_{\Lambda}^TE_2 \\ V_{X}^T \\ V_{X}^TC^T
	\end{bmatrix} &= \mathbbm{b}, \\
	\begin{bmatrix}
	RW_U & HW_X
	\end{bmatrix}
	\begin{bmatrix}
	V_U^T \\ W_X^T
	\end{bmatrix} &= \mathbbm{d}, \\
	\begin{bmatrix}
	W_{\Lambda} & M^TW_{\Lambda} & H^TW_U
	\end{bmatrix}
	\begin{bmatrix}
	V_{\Lambda}^T \\ V_{\Lambda}^TC \\ V_U^T
	\end{bmatrix} &= 0.
	\end{aligned}
	\end{align}
	
	Since using a direct solver would be infeasible, we use an iterative solver, in this case GMRES \cite{Saad1986} to allow for flexibility in choosing a preconditioner, see Section~\ref{sec:Preconditioning}. Algorithm~\ref{alg:LR-GMRES} details a low-rank implementation of GMRES, which leads to low-rank approximate solutions to \eqref{eq:KroneckerWeakSaddle}, making use of \eqref{eq:LRMatrixWeakSaddle}.
	Fundamentally this is the same as a traditional vector-based GMRES with a vector $z$, where instead here we have
	\begin{equation*}
	\mvec{\begin{bmatrix}
		Z_{11}Z_{12}^T \\ Z_{21}Z_{22}^T \\ Z_{31}Z_{32}^T
		\end{bmatrix}} = z.
	\end{equation*}
	Applying the concatenation $X_{k1} = [Y_{k1},\quad Z_{k1}],\; X_{k2} = [Y_{k2},\quad Z_{k2}]$ for $k = 1,2,3$ is equivalent to the vector addition $x = y + z$, since $X_{k1}X_{k2}^T = Y_{k1}Y_{k2}^T + Z_{k1}Z_{k2}^T$ and hence
	
	\begin{equation*}
	x = \mvec{\begin{bmatrix}
		X_{11}X_{12}^T \\ X_{21}X_{22}^T \\ X_{31}X_{32}^T
		\end{bmatrix}} = \mvec{\begin{bmatrix}
		Y_{11}Y_{12}^T + Z_{11}Z_{12}^T \\ Y_{21}Y_{22}^T + Z_{21}Z_{22}^T \\ Y_{31}Y_{32}^T + Z_{31}Z_{32}^T
		\end{bmatrix}} = y+z.
	\end{equation*}
	
	Note that here we employ the same notation as in \cite{Stoll2015}, using the brackets $\{\}$ as a concatenation and truncation operation. Furthermore, after applying the matrix multiplication and the preconditioning, we also truncate the resulting matrices. How this truncation could be implemented is also treated in \cite{Stoll2015}, with options including a truncated singular value decomposition, possibly through Matlab's inbuilt \texttt{svds} function, or a skinny QR factorisation. In the numerical results to follow, we use a modification of the Matlab \texttt{svds} function.
	
	In order to compute the inner product $\langle w,v^{(i)} \rangle$ which arises in GMRES when computing the entries of the Hessenberg matrix (see line 11 in Algorithm~\ref{alg:LR-GMRES}), we make use of the relation between the $\textrm{trace}$ and $\textrm{vec}$ operators:
	\begin{equation*}
	\trace(A^TB) = \mvec{A}^T\mvec{B}.
	\end{equation*}
	Since here 
	\begin{equation*}
	\mvec{\begin{bmatrix}
		W_{11}W_{12}^T \\ W_{21}W_{22}^T \\ W_{31}W_{32}^T
		\end{bmatrix}} = w \quad\text{and}\quad 
	\mvec{\begin{bmatrix}
		V^{(i)}_{11}(V^{(i)}_{12})^T \\ V^{(i)}_{21}(V^{(i)}_{22})^T \\ V^{(i)}_{31}(V^{(i)}_{32})^T
		\end{bmatrix}} = v^{(i)},
	\end{equation*}
	we see that we may compute the inner product $\langle w,v^{(i)} \rangle$ as
	\begin{align}
	\begin{aligned}
	\label{eq:TraceProd1}
	\langle w,v^{(i)} \rangle \;=& \trace\left((W_{11}W_{12}^T)^T(V^{(i)}_{11}(V^{(i)}_{12})^T)\right) + \trace\left((W_{21}W_{22}^T)^T(V^{(i)}_{21}(V^{(i)}_{22})^T)\right) \\&+ 
	\trace\left((W_{31}W_{32}^T)^T(V^{(i)}_{31}(V^{(i)}_{32})^T)\right),
	\end{aligned}
	\end{align}
	by considering the submatrices which make up the vectors $w$ and $v^{(i)}$. Importantly however, the matrices formed in \eqref{eq:TraceProd1} do not exploit the low-rank nature of the submatrices, being $N+1 \times N+1$ matrices. Fortunately, using the properties of the trace operator, we may consider instead:
	\begin{align}
	\label{eq:TraceProd2}
	\langle w,v^{(i)} \rangle \;=&\; \trace\left(W_{11}^TV^{(i)}_{11}(V^{(i)}_{12})^TW_{12}\right) + \trace\left(W_{21}^TV^{(i)}_{21}(V^{(i)}_{22})^TW_{22}\right) \nonumber \\&+ 
	\trace\left(W_{31}^TV^{(i)}_{31}(V^{(i)}_{32})^TW_{32}\right),
	\end{align}
	and hence compute the trace of smaller matrices. In line 11 of Algorithm~\ref{alg:LR-GMRES}, we compute \eqref{eq:TraceProd2} as $\texttt{traceproduct}(W_{11}, W_{12}, W_{21}, W_{22}, W_{31}, W_{32}, V_{11}^{(i)}, V_{12}^{(i)}, V_{21}^{(i)}, V_{22}^{(i)}, V_{31}^{(i)}, V_{32}^{(i)})$.
	
	\begin{algorithm}[!ht]
		\caption{Low-rank GMRES (LR-GMRES)}
		\label{alg:LR-GMRES}
		\begin{algorithmic}
			\STATE Choose $X_{11}^{(0)}, X_{12}^{(0)}, X_{21}^{(0)}, X_{22}^{(0)}, X_{31}^{(0)}, X_{32}^{(0)}.$ 
			\STATE $\{\tilde{X}_{11}, \tilde{X}_{12}, \tilde{X}_{21}, \tilde{X}_{22}, \tilde{X}_{31}, \tilde{X}_{32}\} = \texttt{Amult}(X_{11}^{(0)}, X_{12}^{(0)}, X_{21}^{(0)}, X_{22}^{(0)}, X_{31}^{(0)}, X_{32}^{(0)}).$ 
			\STATE $V_{11} = \{B_{11},\quad -\tilde{X}_{11}\}$, \qquad $V_{12} = \{B_{12},\quad \tilde{X}_{12}\}$,	
			\STATE $V_{21} = \{B_{21},\quad -\tilde{X}_{21}\}$, \qquad $V_{22} = \{B_{22},\quad \tilde{X}_{22}\}$,	
			\STATE $V_{31} = \{B_{31},\quad -\tilde{X}_{31}\}$, \qquad $V_{32} = \{B_{22},\quad \tilde{X}_{32}\}$.	
			\STATE $\xi = [\xi_1, 0, \ldots, 0]$, $\xi_1 =  \sqrt{\texttt{traceproduct}(V_{11}^{(1)}, \ldots,V_{11}^{(1)}, \ldots)}$. 
			\FOR{$k = 1, \ldots $} 
			\STATE  $\{Z_{11}^{(k)}, Z_{12}^{(k)}, Z_{21}^{(k)}, Z_{22}^{(k)}, Z_{31}^{(k)}, Z_{32}^{(k)}\} = \texttt{Aprec}(V_{11}^{(k)}, V_{12}^{(k)}, V_{21}^{(k)}, V_{22}^{(k)}, V_{31}^{(k)}, V_{32}^{(k)})$ 
			\STATE  $\{W_{11}, W_{12}, W_{21}, W_{22}, W_{31}, W_{32}\} = \texttt{Amult}(Z_{11}^{(k)}, Z_{12}^{(k)}, Z_{21}^{(k)}, Z_{22}^{(k)}, Z_{31}^{(k)}, Z_{32}^{(k)}).$ 
			\FOR{$i = 1, \ldots , k$}
			\STATE  $h_{i,k} = \texttt{traceproduct}(W_{11}, \ldots,V_{11}^{(i)}, \ldots)$, 
			\STATE  $W_{11} = \{W_{11},\quad h_{i,k}V_{11}^{(i)}\}$, \qquad $W_{12} = \{W_{12},\quad V_{12}^{(i)}\}$,	
			\STATE  $W_{21} = \{W_{21},\quad h_{i,k}V_{21}^{(i)}\}$, \qquad $W_{22} = \{W_{22},\quad V_{22}^{(i)}\}$,	
			\STATE  $W_{31} = \{W_{31},\quad h_{i,k}V_{31}^{(i)}\}$, \qquad $W_{32} = \{W_{22},\quad V_{32}^{(i)}\}$.	
			\ENDFOR
			\STATE  $h_{k+1,k} = \sqrt{\texttt{traceproduct}(W_{11}, \ldots,W_{11}, \ldots)}$ 
			\STATE  $V_{11}^{(k+1)} = W_{11}/h_{k+1,k}, \qquad V_{12}^{(k+1)} = W_{12}$, 
			\STATE  $V_{21}^{(k+1)} = W_{21}/h_{k+1,k}, \qquad V_{22}^{(k+1)} = W_{22}$,
			\STATE  $V_{31}^{(k+1)} = W_{31}/h_{k+1,k}, \qquad V_{32}^{(k+1)} = W_{32}$.
			\STATE  Apply Givens rotations to $k$th column of $h$, i.e. 
			\FOR{$j = 1,\ldots k-1$}
			\STATE  $\begin{bmatrix}
			h_{j,k} \\ h_{j+1,k}
			\end{bmatrix} = \begin{bmatrix}
			c_j&s_j \\-\bar{s}_j&c_j
			\end{bmatrix}\begin{bmatrix}
			h_{j,k} \\ h_{j+1,k}
			\end{bmatrix}$
			\ENDFOR 
			\STATE  Compute $k$th rotation, and apply to $\xi$ and last column of $h$.
			\STATE  $\begin{bmatrix}
			\xi_{k} \\ \xi_{k+1}
			\end{bmatrix} = \begin{bmatrix}
			c_k&s_k \\ -\bar{s}_k&c_k
			\end{bmatrix}\begin{bmatrix}
			\xi_{k} \\ 0
			\end{bmatrix}$, \qquad 
			$\begin{array}{l}
			h_{k,k} = c_k h_{k,k} + s_k h_{k+1,k}, \\
			h_{k+1,k} = 0.
			\end{array}$
			\IF{$|\xi_{k+1}|$ sufficiently small}
			\STATE  Solve $\tilde{H}\tilde{y} = \xi$, where the entries of $\tilde{H}$ are $h_{i,k}$. 
			\STATE  $Y_{11} = \{\tilde{y}_{1} V_{11}^{(1)}, \;\ldots\;, \;\; \tilde{y}_{k} V_{11}^{(k)}\}, \qquad Y_{12} = \{\tilde{y}_{1} V_{12}^{(1)}, \;\ldots\;, \;\; \tilde{y}_{k} V_{12}^{(k)}\}$  
			\STATE  $Y_{21} = \{\tilde{y}_{1} V_{11}^{(1)}, \;\ldots\;, \;\; \tilde{y}_{k} V_{21}^{(k)}\}, \qquad Y_{22} = \{\tilde{y}_{1} V_{22}^{(1)}, \;\ldots\;, \;\; \tilde{y}_{k} V_{22}^{(k)}\}$ 
			\STATE  $Y_{31} = \{\tilde{y}_{1} V_{31}^{(1)}, \;\ldots\;, \;\; \tilde{y}_{k} V_{31}^{(k)}\}, \qquad Y_{32} = \{\tilde{y}_{1} V_{32}^{(1)}, \;\ldots\;, \;\; \tilde{y}_{k} V_{32}^{(k)}\}$ 
			\STATE  $\{\tilde{Y}_{11}, \tilde{Y}_{12}, \tilde{Y}_{21}, \tilde{Y}_{22}, \tilde{Y}_{31}, \tilde{Y}_{32}\} = \texttt{Aprec}(Y_{11}, Y_{12}, Y_{21}, Y_{22}, Y_{31}, Y_{32})$ 
			\STATE  $X_{11} = \{X_{11}^{(0)}, \;\; \tilde{Y}_{11}\}, \qquad X_{12} = \{X_{12}^{(0)}, \;\; \tilde{Y}_{12}\}$  
			\STATE  $X_{21} = \{X_{21}^{(0)}, \;\; \tilde{Y}_{21}\}, \qquad X_{22} = \{X_{22}^{(0)}, \;\; \tilde{Y}_{22}\}$ 
			\STATE  $X_{31} = \{X_{31}^{(0)}, \;\; \tilde{Y}_{31}\}, \qquad X_{32} = \{X_{32}^{(0)}, \;\; \tilde{Y}_{32}\}$ 
			\STATE \textbf{break}
			\ENDIF
			\ENDFOR
		\end{algorithmic}
	\end{algorithm}
	
	The matrix vector multiplication $Az$ in traditional GMRES, is implemented in LR-GMRES by considering the low-rank form of the saddle point equations generated in \eqref{eq:LRMatrixWeakSaddle}. The concatenation is explicitly written in Algorithm~\ref{alg:Amult} and is denoted \texttt{Amult} in Algorithm~\ref{alg:LR-GMRES}.
	
	\begin{algorithm}[!ht]
		\caption{Matrix multiplication (\texttt{Amult})}
		\label{alg:Amult}
		\begin{algorithmic}
			\REQUIRE $W_{11}, W_{12}, W_{21}, W_{22}, W_{31}, W_{32}$ 
			\ENSURE $Z_{11}, Z_{12}, Z_{21}, Z_{22}, Z_{31}, Z_{32}$ 
			\STATE 
			\begin{tabular}{ll}
				$Z_{11} = [BW_{11},\quad QW_{11},\quad W_{31},\quad MW_{31}]$, 
				&$Z_{12} = [E_1W_{12},\quad E_2W_{12},\quad W_{32},\quad CW_{32}]$, \\
				$Z_{21} = [RW_{21}, \quad HW_{31}]$, 
				&$Z_{21} = [W_{22}, \quad W_{32}]$, \\
				$Z_{31} = [W_{11}, \quad M^TW_{11}, \quad H^TW_{21}]$, 
				&$Z_{32} = [W_{12},\quad C^TW_{12}, \quad W_{22}]$
			\end{tabular}		
		\end{algorithmic}
	\end{algorithm}
	
	Note that we have considered traditional GMRES when implementing LR-GMRES, however it would require only a small modification to allow for restarted GMRES.
	All that remains to consider is preconditioning LR-GMRES, which is implemented in Algorithm~\ref{alg:LR-GMRES} through the \texttt{Aprec} function.
	
	\subsection{Preconditioning LR-GMRES}
	\label{sec:Preconditioning}
	
	We return to the saddle point problem in \eqref{eq:WeakSaddle}. Many approaches exist for preconditioning saddle point problems, a number of which are detailed in \cite{Benzi2005,Benzi2008}. However, the data assimilation setting introduces an unusual situation where the $(1,2)$ block of the saddle point matrix is more computationally expensive than the $(1,1)$ block. In \cite{Fisher2011a,Fisher2016} it is noted that the inexact constraint preconditioner \cite{Bergamaschi2007, Bergamaschi2009, Bergamaschi2011} is an effective choice:
	\begin{equation}
	\label{eq:IC_L0}
	\mathcal{P} = \begin{bmatrix}
	D & 0 & \tilde{L} \\ 0 & \mathcal{R} & 0 \\ \tilde{L}^T & 0 & 0
	\end{bmatrix},
	\end{equation}
	provided a good approximation $\tilde{L}$ to $L =I_{N+1} \kron I_n + C \kron M$ is chosen. Using an inexact constraint preconditioner requires the use of GMRES since the resulting system is nonsymmetric.
	
	Two further requirements must be considered when implementing a preconditioner for LR-GMRES. 
	In order to maintain the low-rank structure we wish to write this in Kronecker form, however we must also consider the inverse of the preconditioner. It is the implementation of the inverse in Kronecker form which allows us to write this as a simple matrix multiplication as in \eqref{eq:LRMatrixWeakSaddle} for the saddle point matrix.
	
	We present here a number of different choices of preconditioner for LR-GMRES.
	
	\subsubsection{Inexact Constraint Preconditioner}
	As mentioned above, the inexact constraint preconditioner \cite{Bergamaschi2011} has been seen to be an effective preconditioner for the saddle point formulation of weak constraint 4D-Var \cite{Fisher2011a}, provided a suitable choice of approximation of $L$ is taken.
	
	The inverse of the inexact constraint preconditioner \eqref{eq:IC_L0} is given by 
	\begin{equation}
	\label{eq:IC_L0inv}
	\mathcal{P}^{-1} = \begin{bmatrix}
	0 & 0 & \tilde{L}^{-T} \\ 0 & \mathcal{R}^{-1} & 0 \\ \tilde{L}^{-1} & 0 & -\tilde{L}^{-1}D\tilde{L}^{-T}
	\end{bmatrix},
	\end{equation}
	which includes the term $\tilde{L}^{-1}$. In order to implement this in LR-GMRES, we write $\tilde{L}^{-1}$ in Kronecker form. This restricts the choice of $\tilde{L}$, however taking an approximation $\tilde{L}$ of the form $I_{N+1} \kron I_n + C \kron \tilde{M}$, where $\tilde{M}$ is an approximation to $M$, the structure of $L$ is maintained. Additionally, we can write the inverse in Kronecker form as
	\begin{align}
	\label{eq:Linv}
	\tilde{L}^{-1} &= I_{N+1} \kron I_n - C \kron \tilde{M} + C^2 \kron \tilde{M}^2 - \ldots + C^N \kron \tilde{M}^N \nonumber\\&= I_{N+1} \kron I_n + \sum_{k = 1}^{N} (-1)^{k} C^{k} \kron \tilde{M}^{k}.
	\end{align}
	Despite being able to write this in Kronecker form, this results in an unfeasible number of terms for large $N$, futhermore for close approximations $\tilde{M}$ to the model matrix $M$, the computations are expensive. A possibility is therefore to approximate $\tilde{L}^{-1}$ by truncating \eqref{eq:Linv} after a few terms.
	
	Truncating after one term we obtain the approximation $\tilde{L}^{-1}=I_{n(N+1)}$. Hence in Kronecker form we can then write the resulting inverse of the preconditioner as:
	\begin{equation}
	\label{eq:P_IC_I0}
	\mathcal{P}_I^{-1} = \begin{bmatrix}
	0 & 0 & I_{N+1} \kron I_n \\ 0 & I_{N+1} \kron R^{-1} & 0 \\ I_{N+1} \kron I_n & 0 & -E_1 \kron B^{-1} + E_2 \kron Q^{-1}
	\end{bmatrix}.
	\end{equation}
	
	To illustrate a possible choice of the \texttt{Aprec} function, we present the application of \eqref{eq:P_IC_I0} as Algorithm~\ref{alg:ICprec}.
	\begin{algorithm}[!ht]
		\caption{Inexact constraint preconditioner $\tilde{L}^{-1}=I_{n(N+1)}$ (\texttt{Aprec})}
		\label{alg:ICprec}
		\begin{algorithmic}
			\REQUIRE $W_{11}, W_{12}, W_{21}, W_{22}, W_{31}, W_{32}$ 
			\ENSURE $Z_{11}, Z_{12}, Z_{21}, Z_{22}, Z_{31}, Z_{32}$ 
			\STATE 
			\begin{tabular}{ll}
				$Z_{11} = W_{31}$, & $Z_{12} = W_{32}$, \\
				$Z_{21} = R^{-1}W_{21}$, & $Z_{21} = W_{22}$, \\
				$Z_{31} = [W_{11}, \quad -B^{-1}W_{31}, \quad -Q^{-1}W_{31}]$, & $Z_{32} = [W_{12},\quad E_1W_{32}, \quad E_2W_{32}]$ 
			\end{tabular}
		\end{algorithmic}
	\end{algorithm}
	
	If we take $\tilde{M} = I_n$ we may consider the approximation $\hat{L} = I_{N+1} \kron I_n + C \kron I_n$. Truncating the resulting inverse after two terms we compute that the Kronecker inverse of the preconditioner is
	\begin{equation}
	\label{eq:P_IC_L0}
	\hat{\mathcal{P}}_{\hat{L}}^{-1} = \begin{bmatrix}
	0 & 0 & I \kron I - C \kron I \\
	0 & I \kron R^{-1} & 0 \\
	I \kron I - C \kron I & 0 & J
	\end{bmatrix},
	\end{equation}
	where $J = - (I \kron I - C \kron I)(E_1 \kron B^{-1})(I \kron I - C^T \kron I) - (I \kron I - C \kron I)(E_2 \kron Q^{-1})(I \kron I - C^T \kron I)$, and we drop the subscripts for the identities.
	
	An alternative approach is to consider an inexact constraint preconditioner where we approximate $\mathcal{H}$ in \eqref{eq:WeakSaddle} in addition to $L$. In this example we approximate $L$ by $\tilde{L} = I$, and using the exact $\mathcal{H}$, we obtain
	\begin{equation}
	\label{eq:IC_IH}
	\mathcal{P}_{I\mathcal{H}} = \begin{bmatrix}
	D & 0 & I \\ 0 & \mathcal{R} & \mathcal{H} \\ I & \mathcal{H}^T & 0
	\end{bmatrix}.
	\end{equation}
	The inverse of which is
	\begin{equation}
	\label{eq:IC_IHinv}
	\mathcal{P}_{I\mathcal{H}}^{-1} = \begin{bmatrix}
	\mathcal{H}^T\mathcal{F}\mathcal{H} & -\mathcal{H}^T\mathcal{F} & I - \mathcal{H}^T\mathcal{F}\mathcal{H}D\\ 
	-\mathcal{F}\mathcal{H} & \mathcal{F} & \mathcal{F}\mathcal{H}D \\ 
	I - D\mathcal{H}^T\mathcal{F}\mathcal{H} & D\mathcal{H}^T\mathcal{F} & D\mathcal{H}^T\mathcal{F}\mathcal{H}D - D
	\end{bmatrix},
	\end{equation}
	where $\mathcal{F} = (\mathcal{H}D\mathcal{H}^T + \mathcal{R})^{-1} = (E_1\kron (HBH^T + R)^{-1}) + (E_2\kron (HQH^T + R)^{-1})$. 
	If $H$ is computationally expensive (such as if $H$ is not a simple interpolatory observation operator), this choice of preconditioner may prove unfeasible.
	
	\subsubsection{Schur Complement Preconditioners}
	An alternative choice of preconditioner is a Schur complement preconditioner, such as the block diagonal preconditioner
	\begin{equation}
	\label{eq:SCpreconditioner}
	\mathcal{P}_D = \begin{bmatrix}
	D & 0 & 0 \\ 0 & \mathcal{R} & 0 \\0 & 0 & \tilde{\mathcal{S}}
	\end{bmatrix},
	\end{equation}
	where $\tilde{\mathcal{S}}$ is an approximation to the Schur-complement
	\begin{equation*}
	\mathcal{S} = -L^TD^{-1}L - \mathcal{H}^TR^{-1}\mathcal{H}.
	\end{equation*}
	This choice of preconditioner is used in \cite{Stoll2015}, and allows the use of LR-MINRES, though in Section~\ref{sec:PrecComp} we use LR-GMRES to compare the different choices as in the full-rank case, GMRES and MINRES are theoretically equivalent for symmetric systems.
	
	As an approximation to the Schur complement we consider
	\begin{equation}
	\label{eq:SchurApprox}
	\tilde{\mathcal{S}} = -\tilde{L}^TD^{-1}\tilde{L},
	\end{equation}
	the inverse of which, $\tilde{\mathcal{S}}^{-1} = -\tilde{L}^{-1}D\tilde{L}^{-T}$ is familiar as the $(3,3)$ term in the inexact constraint preconditioner inverse \eqref{eq:IC_L0inv}. As such we must approximate this by truncating the expansion of $\tilde{L}^{-1}$ \eqref{eq:Linv} as before. Considering the approximation $\hat{L} = I_{N+1} \kron I_n + C \kron I_n$ and truncating after two terms as before, the block diagonal Schur complement preconditioner may be implemented in the same way as the inexact constraint preconditioner \eqref{eq:P_IC_L0} above. This results in 
	\begin{equation}
	\label{eq:P_SC_L0}
	\mathcal{P}_{D\hat{L}}^{-1} = \begin{bmatrix}
	E_1 \kron B^{-1} + E_2 \kron Q^{-1} & 0 & 0\\
	0 & I \kron R^{-1} & 0 \\
	0 & 0 & J
	\end{bmatrix},
	\end{equation}
	where $J = - (I \kron I - C \kron I)(E_1 \kron B^{-1})(I \kron I - C^T \kron I) - (I \kron I - C \kron I)(E_2 \kron Q^{-1})(I \kron I - C^T \kron I)$ as before.

	An alternative method for implementing the Schur complement approximation \eqref{eq:SchurApprox} in a low-rank form is detailed in \cite{Stoll2015}. Instead of truncating the resulting inverse, and applying the technique used in Algorithm~\ref{alg:ICprec}, the relationship between the Kronecker product and Sylvester equations is exploited.
	In order to solve $\tilde{\mathcal{S}}Z_{31}Z_{32}^T = W_{31}W_{32}^T$, the Kronecker form
	\begin{equation*}
	-(I \kron I + C^T \kron \tilde{M}^T)(E_1 \kron B^{-1} + E_2 \kron Q^{-1})(I \kron I + C \kron \tilde{M}) \mvec{Z_{31}Z_{32}^T} = \mvec{W_{31}W_{32}^T},
	\end{equation*}
	is written as two consecutive Sylvester equations. These resulting Sylvester equations are solved one after the other using a low-rank solver such as an ADI \cite{Benner2009, Benner2014} or Krylov \cite{Simoncini2007} method to generate a low-rank approximation $X_{31}X_{32}^T$. It is this approach which we employ in our numerical implementations in Section~\ref{sec:PrecComp}.
	
	An alternative Schur complement preconditioner is the block triangular Schur complement preconditioner, which requires the use of LR-GMRES unlike the block diagonal one above. This choice uses approximations to $L$, $\mathcal{H}$, and the Schur complement $\mathcal{S}$,
	\begin{equation}
	\label{eq:SCpreconditionerT}
	\mathcal{P}_T = \begin{bmatrix}
	D & 0 & \tilde{L} \\ 0 & \mathcal{R} & \tilde{\mathcal{H}} \\0 & 0 & \tilde{\mathcal{S}}
	\end{bmatrix}.
	\end{equation}
	
	When inverted, unlike the other preconditioners we have considered, this maintains a term containing $\tilde{L}$, in addition to the $\tilde{L}^{-1}$ in the Schur complement approximation inverse. Taking the same approximation to $\mathcal{S}$ as above, we obtain the inverse
	\begin{equation}
	\label{eq:SCpreconditionerTinv}
	\mathcal{P}_T^{-1} = \begin{bmatrix}
	D^{-1} & 0 & -D^{-1}\tilde{L}\tilde{\mathcal{S}}^{-1} \\ 0 & \mathcal{R}^{-1} & -\mathcal{R}^{-1}\tilde{\mathcal{H}}\tilde{\mathcal{S}}^{-1} \\0 & 0 & \tilde{\mathcal{S}}^{-1}
	\end{bmatrix}.
	\end{equation}
	In order to implement this preconditioner, \eqref{eq:SCpreconditionerTinv} must be described in Kronecker form, approximating $\tilde{\mathcal{S}}^{-1}$ by truncation or as we use in Section~\ref{sec:PrecComp}, the Sylvester equation approach above.

	\subsubsection{Analysis of preconditioners}
	As mentioned above, whilst there has been investigation into preconditioning saddle point problems such as \cite{Benzi2005,Benzi2008,Bergamaschi2007}, most of these choices assume that the $(1,1)$ block is the computationally expensive one.
	
	Schur complement preconditioners such as the block diagonal and block triangular examples we consider here are detailed in \cite{Benzi2005,Benzi2008}. Using exact matrices for the approximations $\tilde{\mathcal{S}}$, $\tilde{L}$ and $\tilde{\mathcal{H}}$, in \eqref{eq:SCpreconditioner} and \eqref{eq:SCpreconditionerT} results in the preconditioned system having two or three eigenvalues; therefore methods such as MINRES or GMRES converge in at most three steps. However in general, we must consider approximations which reduces the efficacy of the preconditioner. Furthermore, for the data assimilation saddle point problem, these are not necessarily the most appropriate from a computational point of view.
	
	The use of the inexact constraint preconditioner \cite{Bergamaschi2007} in the data assimilation setting is considered in \cite{Fisher2016, Fisher2017, Fisher2011a}, and experimentally has proved effective. Here as the covariance matrices are less computationally expensive, the exact $(1,1)$ block is typically used. Thus using the result in \cite{Bergamaschi2007}, the eigenvalues $\tau$ of the matrix
	\begin{equation}
	\begin{bmatrix}
	D & 0 & \tilde{L} \\ 0 & \mathcal{R} & \tilde{H} \\ \tilde{L}^T & \tilde{H}^T & 0
	\end{bmatrix}^{-1} \begin{bmatrix}
	D & 0 & L \\
	0 & \mathcal{R} & \mathcal{H} \\
	L^T & \mathcal{H}^T & 0
	\end{bmatrix}
	\end{equation}
	are either one (with multiplicity at least $(N+1)(2n+p) - 2\;\mathrm{rank}([L^T,\quad\mathcal{H}^T] - [\tilde{L}^T,\quad \tilde{\mathcal{H}}^T]))$ or bounded by
	\begin{equation*}
	|\tau - 1| \le \frac{\|[L^T,\quad \mathcal{H}^T] - [\tilde{L}^T,\quad \tilde{\mathcal{H}}^T]\|}{\tilde{\sigma}_1},
	\end{equation*}
	where $\tilde{\sigma}_1$ is the smallest singular value of $[\tilde{L}^T,\quad \tilde{\mathcal{H}}^T]$.

	When considering the exact approximation $\tilde{L} = L$, and taking $\tilde{\mathcal{H}} = 0$, the resulting preconditioned system has eigenvalues
	\begin{equation*}
	\tau = 1 \pm \sqrt{\frac{\mu^T\mathcal{H}L^{-1}DL^{-T}\mathcal{H}^T\mu}{\mu^T\mathcal{R}\mu}} \;i
	\end{equation*}
	where $\mu \in \mathbb{R}^{(N+1)p}$. Using the properties of the Rayleigh quotient, we know that the eigenvalues are on a line parallel to the imaginary axis through 1, where the maximum distance from the real axis is given by \begin{equation*}
	\sqrt{\frac{\lambda_{\mathrm{max}}(\mathcal{H}L^{-1}DL^{-T}\mathcal{H}^T)}{\lambda_{\mathrm{min}}(\mathcal{R})}}.
	\end{equation*}
	
	Experimental results in \cite{Fisher2011a} demonstrate that when an approximation is taken for $\tilde{L}$, the eigenvalues are clustered in a cloud surrounding $\tau = 1$ with the size of this cloud likely depending on the accuracy of the chosen approximation.
	
	\section{Numerical Results}
	\label{sec:Numerical}
	In this section we present numerical results using LR-GMRES. (For preconditioning strategies see Section~\ref{sec:PrecComp}). We use 20 iterations of LR-GMRES with a tolerance of $10^{-6}$. During the algorithm where we truncate the matrices after concatenation and applying \texttt{Amult}, we use a truncation tolerance of $10^{-8}$. We present examples with different choices of reduced rank $r$.
	
	\subsection{One-dimensional advection-diffusion system}
	\label{sec:1DAD}
	As a first example, let us consider the one-dimensional (linear) advection-diffusion problem, defined as:
	\begin{equation}
	\label{eq:1DAD}
	\frac{\partial}{\partial t}u(x,t) = c_{d} \frac{\partial^2}{\partial x^2} u(x,t) + c_{a} \frac{\partial}{\partial x}u(x,t)
	\end{equation}
	for $x \in [0,1]$, $t \in (0,T)$, subject to the boundary and initial conditions
	\begin{align*}
	u(0,t) &= 0, &t \in (0,T) \\
	u(1,t) &= 0, &t \in (0,T) \\
	u(x,0) &= u_0(x), &x \in [0,1].
	\end{align*}
	We solve this system with a centered difference scheme for $u_x$ and $u_t$, and a Crank-Nicolson scheme \cite{Crank1947} for $u_{xx}$, discretising $x$ uniformly with $n=100$, and taking timesteps of size $\Delta t = 10^{-3}$.
	For this example, we set the underlying system to have $c_d = 0.1$, $c_a = 1.4$ and $u_0(x) = \sin(\pi x)$. 
	
	We now consider this example as a data assimilation problem, and compare the solutions obtained both by the saddle point formulation \eqref{eq:WeakSaddle}, and the low-rank approximation using LR-GMRES. We take an assimilation window of $200$ timesteps (giving $N = 199$), followed by a forecast of $800$ timesteps. 
	Thus the resulting linear system \eqref{eq:WeakSaddle} we solve here is of size $(40,000 + 200p)$, where $p$ is the number of observations we take at each timestep. Independent of $p$, the full-rank update $\delta x \in\mathbb{R}^{20,000}$. In contrast the low-rank update is $WV^T$, where $W \in \mathbb{R}^{100 \times r}, V \in \mathbb{R}^{200 \times r}$. For $r = 20$, this requires only $30\%$ of the storage.
	
	In the examples to follow, we compare the full- and low-rank solutions to the data assimilation problem with the background estimate.
	
	\paragraph{Perfect observations}
	First let us suppose we have perfect observations of every state in the assimilation window. Hence $p = 100$, and the size of the saddle point system we consider is $60,000$.
	We take as the background estimate $u_0^b$, a perturbed initial condition with background covariance $B = 0.1 I_{100}$, and for this, and the following examples, we consider a model error with covariance $Q = 10^{-4} I_{100}$.
	
	Figure~\ref{fig:AD1D_oN_R0_BI_Nt} shows the state $u(x,t_a)$ and absolute error $|u^*(x,t_a) - u(x,t_a)|$ for the time $t_a$ immediately after assimilation. We consider the three approaches, denoting the true solution by $u^*$. In Figure~\ref{fig:AD1D_oN_R0_BI_RMSE} we consider the root mean squared error of the approaches, presenting the errors in both the assimilation window, and the forecast.
	\begin{figure}[!ht]
		\begin{subfigure}{.48\textwidth}
			\centering
			\resizebox{\textwidth}{!}{
%
%
\definecolor{mycolor1}{rgb}{0.00000,0.44700,0.74100}%
\definecolor{mycolor2}{rgb}{0.92900,0.69400,0.12500}%
\definecolor{mycolor3}{rgb}{0.49400,0.18400,0.55600}%
\definecolor{mycolor4}{rgb}{0.46600,0.67400,0.18800}%
\begin{tikzpicture}

\begin{axis}[%
width=4.359in,
height=3.651in,
at={(0.731in,0.493in)},
scale only axis,
xmin=0,
xmax=100,
xlabel style={font=\color{white!15!black}},
xlabel={x},
ymin=-0.6,
ymax=1,
ylabel style={font=\color{white!15!black}},
ylabel={State},
axis background/.style={fill=white},
legend style={at={(0.657,0.066)}, anchor=south west, legend cell align=left, align=left, draw=white!15!black}
]
\addplot [color=mycolor1, line width=2.0pt]
  table[row sep=crcr]{%
1	0\\
4	0.00616805835302614\\
6	0.0117194504376812\\
8	0.0185997671773492\\
10	0.0269584712470703\\
12	0.0369378941900749\\
14	0.0486685853828135\\
16	0.0622645699054942\\
18	0.0778187040889406\\
20	0.0953983263913898\\
22	0.115041398782168\\
24	0.136753319485706\\
26	0.160504562201794\\
28	0.186229261126115\\
30	0.213824817440198\\
32	0.243152554273053\\
34	0.274039396716688\\
36	0.306280504679236\\
38	0.339642742358848\\
42	0.408682009370906\\
48	0.513688971731284\\
50	0.547867772307455\\
53	0.597327324446226\\
55	0.62867195146822\\
57	0.658390927243133\\
59	0.686222478577761\\
61	0.711918633396081\\
62	0.72389199793453\\
63	0.735244191281069\\
64	0.745947375667086\\
65	0.755974074749858\\
66	0.765296973633667\\
67	0.773888674229241\\
68	0.781721399543102\\
69	0.788766639321622\\
70	0.794994728121438\\
71	0.800374345314211\\
72	0.804871924732453\\
73	0.808450959591326\\
74	0.811071185944144\\
75	0.812687625206109\\
76	0.81324946316461\\
77	0.812698739334309\\
78	0.810968816451251\\
79	0.807982595267532\\
80	0.803650434530638\\
81	0.79786773002678\\
82	0.790512099741463\\
83	0.781440114438155\\
84	0.770483504160268\\
85	0.757444761192687\\
86	0.742092048732161\\
87	0.724153311750186\\
88	0.703309472110362\\
89	0.679186573730306\\
90	0.651346725238881\\
91	0.619277666940462\\
92	0.582380765699611\\
93	0.539957215324506\\
94	0.491192190852871\\
95	0.435136672503376\\
96	0.37068661860107\\
97	0.296559126143322\\
98	0.211265172451903\\
99	0.11307848114167\\
100	0\\
};
\addlegendentry{True State}

\addplot [color=mycolor2, dashed, line width=2.0pt]
  table[row sep=crcr]{%
1	0.170025278399322\\
2	0.00248547182617642\\
4	0.00854130284282917\\
6	0.0162011274118896\\
8	0.0256466530209423\\
10	0.037044944304796\\
12	0.0505411637811193\\
14	0.0662514080702579\\
16	0.0842560259716834\\
18	0.104593814783314\\
20	0.127257473299096\\
22	0.15219064254029\\
24	0.179286789535539\\
26	0.208390089124805\\
28	0.239298340150029\\
30	0.271767823998957\\
31	0.2885021030053\\
33	0.322782188098032\\
35	0.357879057762119\\
43	0.499843486274145\\
46	0.550903097109384\\
48	0.583523696154458\\
50	0.614741197584593\\
52	0.644359869703621\\
54	0.672212553537392\\
56	0.698157445832592\\
58	0.722073398273849\\
60	0.743853990925487\\
62	0.763400604702966\\
63	0.772305583407856\\
64	0.78061462694734\\
65	0.78831398670296\\
66	0.795388771865177\\
67	0.801822610553714\\
68	0.807597280332004\\
69	0.812692300556364\\
70	0.817084477265951\\
71	0.820747389397994\\
72	0.823650802980723\\
73	0.825759997592513\\
74	0.827034986746625\\
75	0.827429610932285\\
76	0.826890478771801\\
77	0.825355728088795\\
78	0.822753574568381\\
79	0.819000611058399\\
80	0.813999815335464\\
81	0.807638218254198\\
82	0.799784177513985\\
83	0.790284194703702\\
84	0.778959204695852\\
85	0.765600256716638\\
86	0.749963495360007\\
87	0.731764337264082\\
88	0.710670724929585\\
89	0.686295323009929\\
90	0.658186504093905\\
91	0.625817950255069\\
92	0.588576673149561\\
93	0.545749228856749\\
94	0.49650587358984\\
95	0.439882372424321\\
96	0.374759134819399\\
97	0.299837307407813\\
98	0.213611405707809\\
99	0.114338011411519\\
100	-0.567527280751037\\
};
\addlegendentry{No assimilation}

\addplot [color=mycolor3, line width=2.0pt]
  table[row sep=crcr]{%
1	0\\
2	-0.0150685556121175\\
3	0.0076125518657193\\
4	0.00100668455947073\\
5	0.0171928766169884\\
6	0.0204864344224092\\
7	-0.00314754749007307\\
8	0.0231684564132735\\
9	0.0302815785330353\\
10	0.00762876649211819\\
11	0.029452766906374\\
12	0.0363193603378846\\
13	0.0425834706797161\\
14	0.026194621942409\\
15	0.065910174235853\\
16	0.0702563387265513\\
17	0.0585949180931209\\
18	0.0708607433602566\\
19	0.0780724727877384\\
20	0.0820032314915693\\
21	0.0814106597892135\\
22	0.116349396195048\\
23	0.117397818580073\\
24	0.126635692693228\\
25	0.163950554333823\\
26	0.161304142569875\\
27	0.186782792505952\\
28	0.181492718441092\\
29	0.198132091782696\\
30	0.213097263776405\\
31	0.237519353045116\\
32	0.255991033505211\\
33	0.254433406742379\\
34	0.27631302731416\\
35	0.311071812597987\\
36	0.304402371530713\\
37	0.329666071586843\\
38	0.341054619644012\\
39	0.367431254041875\\
40	0.368554894368785\\
41	0.405339260440527\\
42	0.421590415678892\\
43	0.425501588091365\\
44	0.44493564392026\\
45	0.450340514153012\\
46	0.495708267769558\\
47	0.509486255639445\\
48	0.511264093225591\\
49	0.523368890525319\\
50	0.553302860560962\\
51	0.574263775759249\\
52	0.582668868106694\\
53	0.581622482461952\\
54	0.616019373650275\\
55	0.638688923750223\\
56	0.632323575344756\\
57	0.671857610172097\\
58	0.689591960543993\\
59	0.691535283281681\\
60	0.708121942433209\\
61	0.71409954328243\\
62	0.729946153212012\\
63	0.733528571553492\\
64	0.749941194357177\\
65	0.747383137511449\\
66	0.770296582981999\\
67	0.788945195974861\\
68	0.789774205135629\\
69	0.794598910829976\\
70	0.78862752855639\\
71	0.805435498753695\\
73	0.803316973746519\\
74	0.821110196269402\\
75	0.84180178043998\\
76	0.821249763143342\\
77	0.803434142467154\\
78	0.801347214717254\\
79	0.81775441754435\\
80	0.795613788205245\\
81	0.790298277861908\\
82	0.779202788338523\\
83	0.782719087672518\\
84	0.778454911297544\\
85	0.751562078658637\\
86	0.750095891837034\\
87	0.721195572548481\\
88	0.708292938018644\\
89	0.690649609616884\\
90	0.661312315734094\\
91	0.619646645203304\\
93	0.546919907072038\\
95	0.421600543764697\\
96	0.385618813850613\\
97	0.283502614600266\\
98	0.218010888490369\\
99	0.118626957094705\\
100	0\\
};
\addlegendentry{Full-rank}

\addplot [color=mycolor4, line width=2.0pt]
  table[row sep=crcr]{%
1	0\\
2	-0.0146719663362518\\
3	0.0120346001447302\\
4	0.0055661244128089\\
5	0.0142139488948629\\
6	0.0174222255676\\
7	-0.00343223848364005\\
8	0.0276812437285088\\
9	0.0317757921558268\\
10	0.0117201263987567\\
11	0.0324180427237337\\
12	0.0374603588855962\\
13	0.041637812853395\\
14	0.0298629545360143\\
15	0.0633043463979561\\
16	0.0703645688946324\\
17	0.0560776234161722\\
18	0.0746643109481795\\
19	0.0816434312002201\\
20	0.0844493374057009\\
21	0.0861676034374455\\
22	0.116480514287844\\
23	0.115792562233352\\
24	0.123432514994107\\
25	0.165921250430529\\
26	0.157891542288766\\
27	0.187431372381212\\
28	0.184781730783854\\
29	0.198893002962038\\
30	0.212133663003684\\
31	0.240064880205736\\
32	0.252392725311708\\
33	0.254673658952854\\
34	0.270554292381618\\
35	0.306494189884077\\
36	0.297628534029457\\
37	0.33010986898779\\
38	0.338698115710031\\
39	0.359444923343034\\
40	0.363944538200471\\
41	0.398237618632294\\
42	0.417328499092633\\
43	0.422994707987485\\
44	0.443130451750562\\
45	0.449787120820716\\
46	0.496665447780572\\
47	0.5108477416195\\
48	0.51323268938512\\
49	0.523638413789755\\
50	0.553864361983187\\
51	0.578689681269239\\
52	0.580890453952705\\
53	0.580263878781878\\
54	0.613005992274566\\
55	0.640601811375717\\
56	0.62876407106593\\
57	0.674135362657609\\
58	0.68367629740726\\
59	0.68736245678285\\
60	0.702920816535993\\
61	0.709962697218842\\
62	0.723778705350284\\
63	0.736571162181747\\
64	0.751811358395571\\
65	0.747618875945108\\
66	0.77171235410917\\
67	0.787345672457477\\
68	0.783459401377797\\
69	0.793509910724865\\
70	0.78993432254245\\
71	0.806505191826062\\
72	0.807195357598971\\
73	0.807132971845334\\
74	0.820167982795084\\
75	0.843344635192409\\
76	0.818689752103154\\
77	0.805581640374726\\
78	0.805876318852697\\
79	0.816679586311921\\
80	0.801211224462918\\
81	0.795460598889122\\
82	0.782005483216395\\
83	0.781989231682189\\
84	0.774375415106093\\
85	0.750015697253659\\
86	0.757336818252426\\
87	0.722908336521101\\
88	0.712247942236573\\
89	0.6959386251295\\
90	0.66817981017509\\
91	0.620471362451227\\
92	0.582716756023274\\
93	0.547862407368569\\
94	0.489575329442502\\
95	0.423256481503657\\
96	0.387219200603695\\
97	0.278482366201629\\
98	0.221863043561854\\
99	0.12050719026972\\
100	0\\
};
\addlegendentry{Low-rank}

\end{axis}
\end{tikzpicture}%
}
			\caption{State $u(x,t_a)$}
		\end{subfigure}
		\begin{subfigure}{.48\textwidth}
			\centering
			\resizebox{\textwidth}{!}{
%
%
\definecolor{mycolor1}{rgb}{0.92900,0.69400,0.12500}%
\definecolor{mycolor2}{rgb}{0.49400,0.18400,0.55600}%
\definecolor{mycolor3}{rgb}{0.46600,0.67400,0.18800}%
\begin{tikzpicture}

\begin{axis}[%
width=4.359in,
height=3.651in,
at={(0.731in,0.493in)},
scale only axis,
xmin=0,
xmax=100,
xlabel style={font=\color{white!15!black}},
xlabel={x},
ymin=0,
ymax=0.0757034505375486,
ylabel style={font=\color{white!15!black}},
ylabel={Error},
axis background/.style={fill=white},
legend style={legend cell align=left, align=left, draw=white!15!black}
]
\addplot [color=mycolor1, dashed, line width=2.0pt]
  table[row sep=crcr]{%
1.51231232854519	0.0832737955913103\\
2	0.000692084614243527\\
3	0.00148148203700771\\
4	0.00237324448980303\\
5	0.00337196383907212\\
6	0.00448167697419422\\
7	0.00570577045253629\\
8	0.00704688584359303\\
9	0.00850682731122276\\
10	0.0100864730577399\\
11	0.0117856923140067\\
12	0.0136032695910444\\
13	0.0155368379091527\\
14	0.0175828226874444\\
15	0.0197363979082468\\
16	0.0219914560661891\\
17	0.0243405932708072\\
18	0.0267751106943592\\
19	0.0292850333457721\\
20	0.0318591469077063\\
21	0.0344850531028271\\
22	0.0371492437581225\\
26	0.047885526923011\\
27	0.0505071407787199\\
28	0.0530690790239134\\
29	0.0555535940122809\\
30	0.0579430065587587\\
31	0.0602199022943779\\
32	0.0623673309860635\\
33	0.0643690058838899\\
34	0.0662095000640903\\
35	0.067874436698645\\
36	0.069350670201402\\
37	0.0706264552799212\\
38	0.0716916010596691\\
39	0.0725376076425022\\
40	0.0731577827113057\\
41	0.0735473360929433\\
42	0.0737034505375505\\
43	0.0736253273557992\\
44	0.0733142059716698\\
45	0.0727733568851789\\
46	0.0720080479919005\\
47	0.0710254846608649\\
48	0.0698347244231741\\
49	0.0684465675584676\\
50	0.0668734252771372\\
51	0.0651291675740708\\
52	0.0632289531655488\\
53	0.061189044209641\\
54	0.0590266087438067\\
55	0.0567595139483927\\
56	0.0544061134570484\\
57	0.0519850319834774\\
58	0.0495149505171071\\
62	0.039508606768436\\
63	0.0370613921267875\\
64	0.03466725128024\\
65	0.032339911953116\\
66	0.0300917982315099\\
67	0.0279339363244588\\
68	0.0258758807889024\\
69	0.0239256612347276\\
70	0.0220897491445129\\
71	0.020373044083783\\
72	0.0187788782482698\\
73	0.0173090380011871\\
74	0.0159638008024672\\
75	0.0147419857261752\\
76	0.0136410156071918\\
77	0.0126569887544861\\
78	0.0117847581171304\\
79	0.0110180157908673\\
80	0.0103493808048256\\
81	0.00977048822741722\\
82	0.00927207777252193\\
83	0.00884408026553274\\
84	0.008475700535584\\
85	0.00815549552395112\\
86	0.00787144662783135\\
90	0.00683977885502429\\
91	0.00654028331459244\\
92	0.0061959074499498\\
93	0.00579201353222913\\
94	0.00531368273696842\\
95	0.00474569992094587\\
96	0.00407251621832927\\
97	0.00327818126449131\\
98	0.00234623325590633\\
99	0.00125953026983439\\
99.1448330145091	0.0832737955913103\\
};
\addlegendentry{No assimilation}

\addplot [color=mycolor2, line width=2.0pt]
  table[row sep=crcr]{%
1	0\\
2	0.0168619428240646\\
3	0.00376941169068346\\
4	0.00516137379355541\\
5	0.00840574032226016\\
6	0.00876698398472797\\
7	0.0181315852643991\\
8	0.00456868923592424\\
9	0.00769637421423397\\
10	0.0193297047549521\\
11	0.00228433489982649\\
12	0.000618533852190239\\
13	6.70866322138863e-06\\
14	0.0224739634404045\\
15	0.0106831108149947\\
16	0.00799176882105712\\
17	0.0111970953913669\\
18	0.00695796072869825\\
19	0.00827975426223304\\
20	0.0133950948998205\\
21	0.0235501388994379\\
22	0.00130799741288001\\
23	0.00824193993022959\\
24	0.0101176267924785\\
25	0.0155732735547502\\
26	0.000799580368081365\\
27	0.0136570088126717\\
28	0.00473654268502344\\
29	0.00166892686038977\\
30	0.000727553663793401\\
31	0.00923715233420808\\
32	0.0128384792321583\\
33	0.00397977547176254\\
34	0.00227363059745755\\
35	0.0210671915345131\\
36	0.00187813314852292\\
37	0.00682902264881591\\
38	0.00141187728516456\\
39	0.0107665503155658\\
40	0.00531393723885287\\
41	0.0141193017465042\\
42	0.0129084063079858\\
43	0.000716570826980956\\
44	0.00114465116280371\\
45	0.0110221500363821\\
46	0.0168132186520751\\
47	0.0131363593408906\\
48	0.00242487850569262\\
49	0.00750530533088067\\
50	0.00543508825350614\\
51	0.00963146874052256\\
52	0.00153795156863623\\
53	0.0157048419842738\\
54	0.00283342885668958\\
55	0.0100169722820027\\
56	0.0114277570307877\\
57	0.0134666829289785\\
58	0.0170335127872505\\
59	0.00531280470391948\\
60	0.00876948485007745\\
61	0.00218090988633435\\
62	0.00605415527748221\\
63	0.00171561972759093\\
64	0.00399381869009119\\
65	0.00859093723840942\\
66	0.00499960934833155\\
67	0.0150565217456204\\
68	0.00805280559252708\\
69	0.00583227150835341\\
70	0.00636719956504805\\
71	0.00506115343948466\\
72	0.000242492601827848\\
74	0.0100390103252579\\
75	0.0291141552338701\\
76	0.00800029997873253\\
77	0.00926459686715475\\
78	0.00962160173399695\\
79	0.00977182227681794\\
80	0.0080366463253938\\
81	0.00756945216487281\\
82	0.0113093114029255\\
83	0.00127897323434922\\
84	0.00797140713726208\\
85	0.00588268253405033\\
86	0.00800384310485924\\
87	0.00295773920170461\\
88	0.00498346590828191\\
89	0.0114630358865782\\
90	0.00996559049521295\\
91	0.000368978262827113\\
92	0.00118966653397479\\
93	0.00696269174753184\\
94	0.00718150481368696\\
95	0.0135361287386786\\
96	0.0149321952495427\\
97	0.0130565115430556\\
98	0.0067457160384663\\
99	0.00554847595303443\\
100	0\\
};
\addlegendentry{Full-rank}

\addplot [color=mycolor3, line width=2.0pt]
  table[row sep=crcr]{%
1	0\\
2	0.0164653535482415\\
3	0.00819145996965176\\
4	0.00060193394027408\\
5	0.00542681260031941\\
6	0.00570277512991879\\
7	0.018416276257966\\
8	0.00908147655114533\\
9	0.00919058783719606\\
10	0.0152383448483846\\
11	0.000680940917391126\\
12	0.000522464695521307\\
13	0.0009389491631282\\
14	0.0188056308468276\\
15	0.00807728297702681\\
16	0.0080999989890671\\
17	0.0137143900683299\\
18	0.00315439314064747\\
19	0.0047087958496661\\
20	0.010948988985632\\
21	0.0187931952511917\\
22	0.0014391155056046\\
23	0.0098471962768798\\
24	0.0133208044915705\\
25	0.0175439696513706\\
26	0.00261301991311313\\
27	0.0143055886878898\\
28	0.00144753034226142\\
29	0.00090801568107679\\
30	0.0016911544364433\\
31	0.0117826794947717\\
32	0.00924017103865538\\
33	0.00373952326137328\\
34	0.00348510433505567\\
35	0.0164895688205604\\
36	0.00865197064970857\\
37	0.00727282004982044\\
38	0.000944626648802682\\
39	0.00278021961672437\\
40	0.00992429340722367\\
41	0.0070176599382421\\
42	0.00864648972171267\\
43	0.00322345093073295\\
44	0.000660541006880067\\
45	0.011575543368636\\
46	0.0177703986632025\\
47	0.0144978453210598\\
48	0.00045628234605033\\
49	0.00723578206643083\\
50	0.00599658967573191\\
51	0.014057374250541\\
52	0.000240462585324508\\
53	0.0170634456643626\\
54	0.000179952519061999\\
55	0.0119298599074966\\
56	0.014987261309571\\
57	0.0157444354145753\\
58	0.0111178496504181\\
59	0.0011399782050745\\
60	0.0035683589529043\\
61	0.00195593617718259\\
62	0.000113292584231317\\
63	0.00132697090072043\\
64	0.00586398272861288\\
65	0.00835519880465085\\
66	0.00641538047545964\\
67	0.0134569982280937\\
68	0.00173800183465289\\
69	0.00474327140311459\\
70	0.00506040557898757\\
71	0.00613084651180884\\
72	0.00232343286650405\\
73	0.00131798774592085\\
74	0.00909679685095455\\
75	0.0306570099862\\
76	0.00544028893845905\\
77	0.00711709895954016\\
78	0.00509249759846853\\
79	0.00869699104444521\\
80	0.00243921006776304\\
81	0.00240713113758773\\
82	0.00850661652512485\\
83	0.000549117244034392\\
84	0.00389191094583907\\
85	0.00742906393908527\\
86	0.0152447695202653\\
87	0.00124497522917011\\
88	0.00893847012619631\\
89	0.016752051399294\\
90	0.0168330849362093\\
91	0.001193695510878\\
92	0.00033599032371967\\
93	0.00790519204414863\\
94	0.00161686141029804\\
95	0.0118801909997188\\
96	0.0165325820025117\\
97	0.018076759941664\\
98	0.0105978711099226\\
99	0.00742870912802118\\
100	0\\
};
\addlegendentry{Low-rank}

\end{axis}
\end{tikzpicture}%
}
			\caption{Error $|u^*(x,t_a) - u(x,t_a)|$}
		\end{subfigure}
		\caption{State and error for time $t_a$ after the assimilation window for 1D advection-diffusion problem with perfect observations.}
		\label{fig:AD1D_oN_R0_BI_Nt}
	\end{figure}
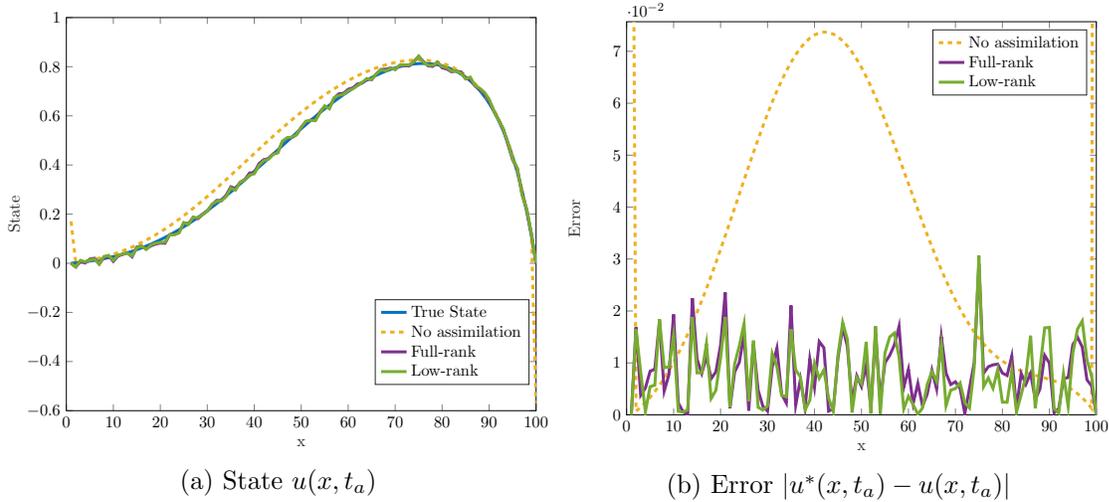
	\begin{figure}[!ht]
		\centering
		\resizebox{.48\textwidth}{!}{\input{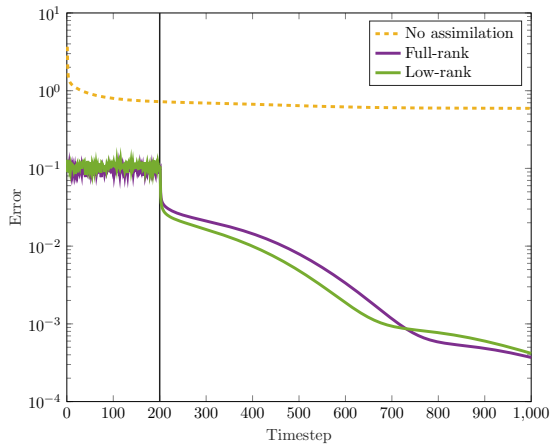}}
		\caption{Root mean squared errors for 1D advection-diffusion data assimilation problem with perfect observations.}
		\label{fig:AD1D_oN_R0_BI_RMSE}
	\end{figure}
	The results show that the low-rank solution matches the full-rank solution very closely, in both the observation window and the forecast. In Figure~\ref{fig:AD1D_oN_R0_BI_Nt}, the low, and full-rank approximations are indistinguishable, with both displaying the same characteristics in the state error plot. Both methods for solving the data assimilation problem result in a superior forecast to the initial guess (without assimilation).
	
	It is worth noting that here the low-rank solution to the data assimilation problem achieves a lower root mean squared error than the full-rank solution for half of the forecast window. Investigating different random seeds, we saw that this was not always the case, though in majority of experiments the two solutions were close. In this example, the full- and low-rank solutions both outperformed the background estimate for all random seeds considered.
	
\begin{figure}[!ht]
		\begin{subfigure}{.48\textwidth}
			\centering
			\resizebox{\textwidth}{!}{
%
%
\definecolor{mycolor1}{rgb}{0.92900,0.69400,0.12500}%
\definecolor{mycolor2}{rgb}{0.49400,0.18400,0.55600}%
\definecolor{mycolor3}{rgb}{0.46600,0.67400,0.18800}%
\begin{tikzpicture}

\begin{axis}[%
width=4.359in,
height=3.651in,
at={(0.731in,0.493in)},
scale only axis,
xmin=0,
xmax=100,
xlabel style={font=\color{white!15!black}},
xlabel={x},
ymin=0,
ymax=0.234861171525923,
ylabel style={font=\color{white!15!black}},
ylabel={Error},
axis background/.style={fill=white},
legend style={legend cell align=left, align=left, draw=white!15!black}
]
\addplot [color=mycolor1, dashed, line width=2.0pt]
  table[row sep=crcr]{%
1	0.114207108891847\\
2	0.000281350604765862\\
4	0.000959797477634083\\
6	0.00179292990434021\\
8	0.00277177161737541\\
11	0.00446080492386614\\
17	0.00800121141843135\\
19	0.00892913721055777\\
20	0.00928326551242264\\
21	0.00954182614891863\\
22	0.00968684098330641\\
23	0.00969950092225247\\
24	0.0095603131162818\\
25	0.00924927052716384\\
26	0.00874604304071624\\
27	0.00803018884887763\\
28	0.00708138435895478\\
29	0.00587967042061166\\
30	0.00440571220387653\\
31	0.00264106962403332\\
32	0.000568474804637731\\
33	0.00182788729227923\\
34	0.00456209924607265\\
35	0.00764623698663058\\
36	0.0110901056188766\\
37	0.0149009899328121\\
38	0.0190834242629876\\
39	0.0236389862618864\\
40	0.0285661189536341\\
41	0.0338599851400971\\
42	0.0395123578411472\\
43	0.045511549969433\\
44	0.0518423858747781\\
45	0.0584862167530673\\
46	0.0654209812106359\\
47	0.0726213115213596\\
48	0.0800586853240901\\
49	0.0877016216997362\\
50	0.0955159197561102\\
51	0.103464937053019\\
53	0.119610273132025\\
55	0.135808390736486\\
56	0.143819619670126\\
57	0.151714045766781\\
58	0.159448192476418\\
59	0.166979260260661\\
60	0.174265539900645\\
61	0.181266809223914\\
62	0.187944706740069\\
63	0.194263075932056\\
64	0.200188274291094\\
65	0.205689441592966\\
66	0.210738722375496\\
67	0.215311438069946\\
68	0.219386204741411\\
69	0.222944992880656\\
70	0.225973126135528\\
71	0.228459216247316\\
72	0.230395031736322\\
73	0.231775298032773\\
74	0.232597426740952\\
75	0.232861171525926\\
76	0.232568207689482\\
77	0.231721631819681\\
78	0.230325376923687\\
79	0.228383537146996\\
80	0.225899594507212\\
81	0.22287553798617\\
82	0.219310862789428\\
83	0.215201434552085\\
84	0.210538199700594\\
85	0.205305719023144\\
86	0.199480496710251\\
87	0.193029071652475\\
88	0.185905831578538\\
89	0.178050503640748\\
90	0.169385267265781\\
91	0.15981142645991\\
92	0.149205569269583\\
93	0.137415131752221\\
94	0.124253272629943\\
95	0.109492952835808\\
96	0.0928601015080659\\
97	0.0740257367849608\\
98	0.0525968961975565\\
99	0.0281062178230229\\
100	0.052538469640325\\
};
\addlegendentry{No assimilation}

\addplot [color=mycolor2, line width=2.0pt]
  table[row sep=crcr]{%
1	0\\
2	0.0215149183640619\\
3	0.0515474960005662\\
4	0.0281822623677783\\
5	0.0129210656638463\\
6	0.0872463928130429\\
7	0.0454786502065048\\
8	0.00304159940399984\\
9	0.0312802235516614\\
10	0.00185311044192815\\
11	0.0878695189757792\\
12	0.0199241589478305\\
13	0.0120111271340306\\
14	0.0209908637763192\\
15	0.042617556715939\\
16	0.122759994988343\\
17	0.0295267726198318\\
18	0.0383147202872038\\
19	0.00187759419164024\\
20	0.0203970993626683\\
21	0.00177704055100492\\
22	0.054655044230941\\
23	0.0263413270958068\\
24	0.0642098293959918\\
25	0.0817601539026924\\
26	0.0269053415280354\\
27	0.0627011747434949\\
28	0.0153325253833998\\
29	0.0173845386476472\\
30	0.0165074088228039\\
31	0.0231885854524876\\
32	0.0127240271883267\\
33	0.0465735213638538\\
34	0.059144838195806\\
35	0.0693745905857526\\
36	0.0554882432501103\\
37	0.0659125207176317\\
38	0.0167790681255724\\
39	0.0450349663954199\\
40	0.00631722485199759\\
41	0.0753347470256784\\
42	0.0133613582853798\\
43	0.0331912646705348\\
44	0.0242171892132177\\
45	0.0432472791813581\\
46	0.0727878373377422\\
47	0.0514078249890133\\
48	0.075889404901119\\
49	0.0205426273278277\\
50	0.0227324179038249\\
51	0.00369251736012188\\
52	0.0383168345275919\\
53	0.0388062008187546\\
54	0.0772076674421243\\
55	0.0715653568801855\\
56	0.0155084231200107\\
57	0.0466630705757325\\
58	0.0491842212170042\\
59	0.0575307886370808\\
60	0.0615960338075041\\
61	0.0684929788936017\\
62	0.0105590714946118\\
63	0.0387871357431777\\
64	0.0111513411161468\\
65	0.0218763106135356\\
66	0.0159220636966779\\
67	0.0179349634251196\\
68	0.0204561119563778\\
69	0.00928491107711693\\
70	0.0127816308273196\\
71	0.0636149221653568\\
72	0.0102351185103942\\
73	0.00692550221339161\\
74	0.0154191179848198\\
75	0.0108294588282405\\
76	0.053259619496103\\
77	0.0361589974542369\\
78	0.0281155580351822\\
79	0.0311462951066517\\
80	0.0440216469769439\\
81	0.105933995891334\\
82	0.10873840047249\\
83	0.0798076100881957\\
84	0.0266588678597799\\
85	0.0189158729853176\\
86	0.0620860832595724\\
87	0.0398647270948089\\
88	0.0124425863554336\\
89	0.0332821503829592\\
90	0.0266020138813872\\
91	0.0715972074474394\\
92	0.042598921154493\\
93	0.0588457825619599\\
94	0.0674598012229097\\
95	0.0341701596344137\\
96	0.118719248632246\\
97	0.0444681253723331\\
98	0.0651398314203817\\
99	0.0279228720307572\\
100	0\\
};
\addlegendentry{Full-rank}

\addplot [color=mycolor3, line width=2.0pt]
  table[row sep=crcr]{%
1	0\\
2	0.00335195101787633\\
3	0.0177500284449224\\
4	0.0332885347053917\\
5	0.0311442267527013\\
6	0.0320856414637802\\
7	0.0564968389534499\\
8	0.0136546102423125\\
9	0.00296661020757938\\
10	0.039976547636229\\
11	0.061061131189831\\
12	0.0471932158447999\\
13	0.00733768856014194\\
14	0.0125844181539634\\
15	0.00353879619616748\\
16	0.0595480494194618\\
17	0.0101506476129316\\
18	0.0104855005004367\\
19	0.0178305605467131\\
20	0.0615686787428587\\
21	0.0238254484876279\\
22	0.0926913151624831\\
23	0.0115699355920214\\
24	0.0225083509095327\\
25	0.0755460959123582\\
26	0.0136510990099055\\
27	0.0695116262010345\\
28	0.0122385253341974\\
29	0.0247160204981896\\
30	0.000679651148658422\\
31	0.0443472113231564\\
32	0.0106716937730482\\
33	0.0564250475265595\\
34	0.048567131306541\\
35	0.0524301868905326\\
36	0.0808599607178451\\
37	0.0355534648914926\\
38	0.00677024913517243\\
39	0.0543100208117977\\
40	0.0621124197596856\\
41	0.0403856638897224\\
42	0.0795176422704458\\
43	0.0636720720571731\\
44	0.0559827803234612\\
45	0.0885652110817148\\
46	0.124307996256675\\
47	0.0642631312962436\\
48	0.0702059030149513\\
49	0.0345768821112955\\
50	0.0648584256859408\\
51	0.0211249700450082\\
52	0.101216144517423\\
53	0.0439114496791007\\
54	0.0531722662097422\\
55	0.0323024730389676\\
56	0.0344308593906675\\
57	0.0151696303520765\\
58	0.0276752680669006\\
59	0.028853648029056\\
60	0.0295003026425604\\
61	0.0869889451649755\\
62	0.00993106466074778\\
63	0.0162333718117793\\
64	0.0376961490582488\\
65	0.0496778700108393\\
66	0.0293979063719547\\
67	0.0667977100787738\\
68	0.0395788779755293\\
69	0.0436085738895855\\
70	0.0461384040497848\\
71	0.0762003333466481\\
72	0.0254047889851989\\
73	0.030482379661251\\
74	0.00411164930098096\\
75	0.0237189674595442\\
76	0.0590848359968703\\
77	0.00370644849981261\\
78	0.03195206387052\\
79	0.037423258495167\\
80	0.0394163429030243\\
81	0.0155589231519997\\
82	0.0148650428199915\\
83	0.00877189448999616\\
84	0.0183781325826544\\
85	0.0550910013505472\\
86	0.00550826542894356\\
87	0.0334280596919143\\
88	0.0191315848459368\\
89	0.0374499967674069\\
90	0.00623495037450539\\
91	0.0431169841023547\\
92	0.0064687108857413\\
93	0.0146756945303821\\
94	0.0421590508644556\\
95	0.000726512055763351\\
96	0.0994447629249748\\
97	0.0296400948706861\\
98	0.0213858886905456\\
99	0.00407339190563505\\
100	0\\
};
\addlegendentry{Low-rank}

\end{axis}
\end{tikzpicture}%
}
			\caption{Error $|u^*(x,t_a) - u(x,t_a)|$}
		\end{subfigure}
		\begin{subfigure}{.48\textwidth}
			\centering
			\resizebox{\textwidth}{!}{\input{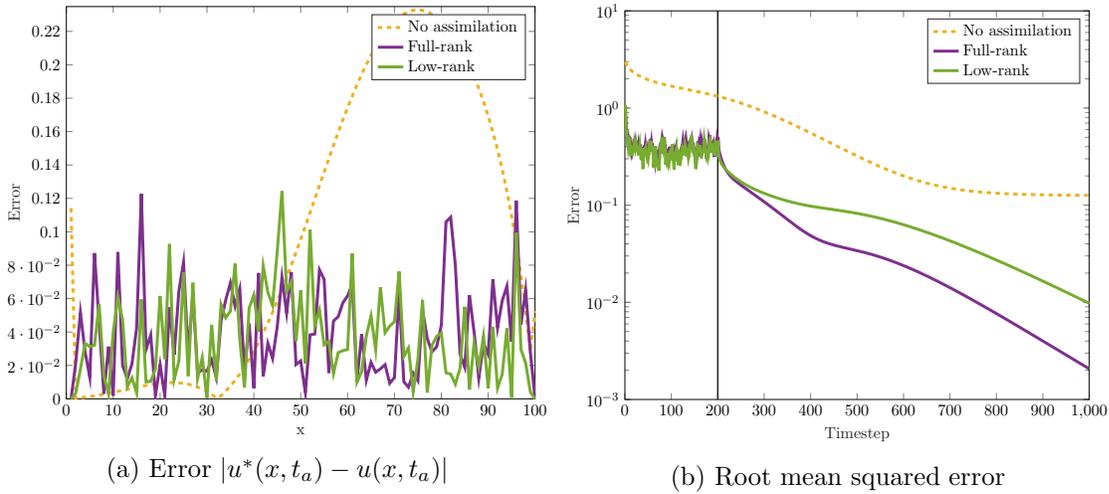}}
			\caption{Root mean squared error}
		\end{subfigure}
		\caption{Error for time $t_a$ after the assimilation window, and root mean squared error for 1D advection-diffusion problem with partial, noisy observations $(r=20)$.}
		\label{fig:AD1D_oN5_RI_Be_E}
	\end{figure}

	\paragraph{Partial, noisy observations}
	Next, we introduce partial noisy observations, taking observations in every fifth component of $u$. These are generated from the truth with covariance $R = 0.01 I_{p}$, for $p=20$, and as such the linear system we consider for this example is of size $44,000$.
	In this example we take for the background error covariance $B_{i,j} = 0.1\exp(\frac{-|i-j|}{50})$, keeping $Q = 10^{-4} I_{100}$ and $r=20$. The resulting errors for three approaches, and the root mean squared errors are shown in Figure~\ref{fig:AD1D_oN5_RI_Be_E}. 
	
	As with the previous example, the state errors of both the full- and low-rank solutions are similar, though here we notice a greater variation between the two than in the previous example.
	Unlike above, when we compare the root mean squared errors of the full- and low-rank approaches, there is a greater disparity between the two, with the full-rank performing significantly better except at the very start of the forecast. Nonetheless the low-rank approximation is superior to using no assimilation.
	
	\paragraph{Different choices of rank}
	We now consider the effect of the chosen rank on the assimilation result. In the previous examples we have considered $r=20$, which resulted in the low-rank approximation to $\delta x$ requiring only $30\%$ of the storage needed for the full-rank solution. Here we consider $r = 5$ (requiring $7.5\%$ of the storage), and $r = 1$ (needing just $1.5\%$), and otherwise keep the setup of the example used in Figure~\ref{fig:AD1D_oN5_RI_Be_E}, with partial, noisy observations unchanged.
	\begin{figure}[!ht]
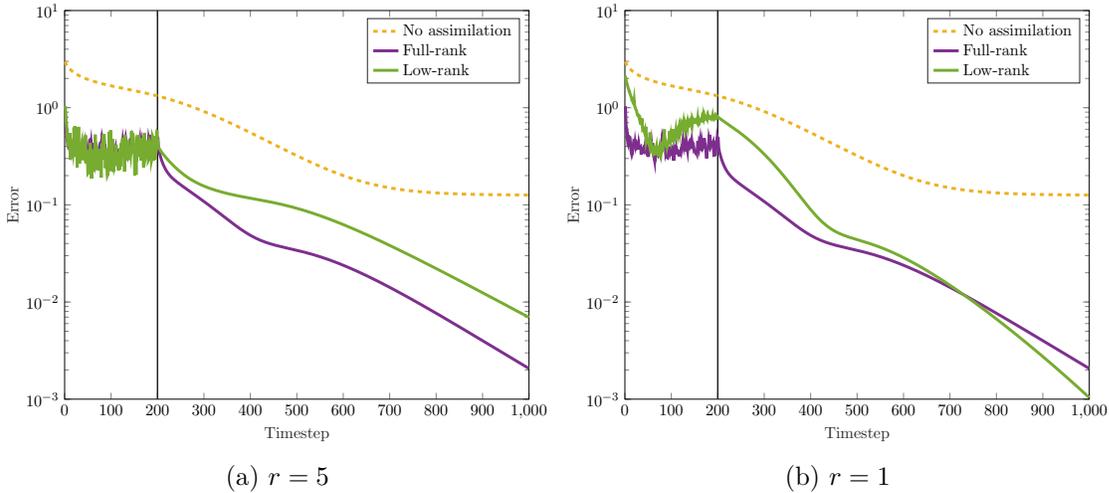

		\begin{subfigure}{.48\textwidth}
			\centering
			\resizebox{\textwidth}{!}{\input{tikzplots/AD1D_oN5_RI_Be_r5_RMSE.tex}}
			\caption{$r=5$}
		\end{subfigure}
		\begin{subfigure}{.48\textwidth}
			\centering
			\resizebox{\textwidth}{!}{\input{tikzplots/AD1D_oN5_RI_Be_r1_RMSE.tex}}
			\caption{$r=1$}
		\end{subfigure}
		\caption{Root mean squared errors for 1D advection-diffusion data assimilation problem with partial, noisy observations.}
		\label{fig:AD1D_oN5_RI_Be_RMSEr}
	\end{figure}
	In Figure~\ref{fig:AD1D_oN5_RI_Be_RMSEr} we obtain a very close forecast from taking $r=5$ to that which we saw from $r=20$, though the assimilation window has greater variation for $r=5$ whilst remaining close to the full-rank solution. In contrast, the behaviour of the root mean squared error for $r=1$ is considerably different to that of the full-rank solution. Despite this, the forecasts for both $r=5$ and $r=1$ are close to the full-rank solution and are comfortably more accurate than using no assimilation. The closeness to the full-rank may be caused by the smoothing properties of this model operator, and the particular random seed, as noted above. Though taking different random seeds results in similar behaviour in majority of cases.
	
	Table~\ref{table:ADcomp} presents the storage requirements for the examples considered in this section. As Figures~\ref{fig:AD1D_oN_R0_BI_Nt}-~\ref{fig:AD1D_oN5_RI_Be_RMSEr} demonstrate, despite the large reduction in the necessary storage for the low-rank approach, it results in close approximations to the full-rank method.
	\begin{table}[!ht]
		\centering
		\begin{tabular}{|l|l|l|l|l|l|l|}
			\hline
			&&&& \multicolumn{2}{|l|}{\# of matrix elements in solution} & \\\cline{5-6}
			n & N & p & rank &full-rank & low-rank & storage reduction\\
			\hline
			100 & 199 & 100 & 20 & 20,000 & 6,000 & 70\% \\
			100 & 199 & 20  & 20 & 20,000 & 6,000 & 70\% \\
			100 & 199 & 20  & 5  & 20,000 & 1,500 & 92.5\% \\		
			100 & 199 & 20  & 1  & 20,000 &   300 & 98.5\% \\	
			\hline
		\end{tabular}
		\caption{Storage requirements for full- and low-rank methods in the 1D advection-diffusion equation examples.}
		\label{table:ADcomp}
	\end{table}
	
	\FloatBarrier
	\subsection{Comparison of preconditioners}
	\label{sec:PrecComp}
	We present here a comparison between different choices of preconditioner for the 1D advection -diffusion equation system in Section~\ref{sec:1DAD}. We consider a small example taking $n=10$, $N = 19$, $p = 4$ with $B_{i,j} = 0.1\exp(\frac{-|i-j|}{50})$, $Q = 10^{-4}I_{10}$, $R = 0.01I_4$. The resulting saddle point matrix is $440 \times 440$. In all the following cases a reduced rank size of $r=5$ is considered, though similar results are obtained when we vary this choice.
	
\begin{figure}[!ht]
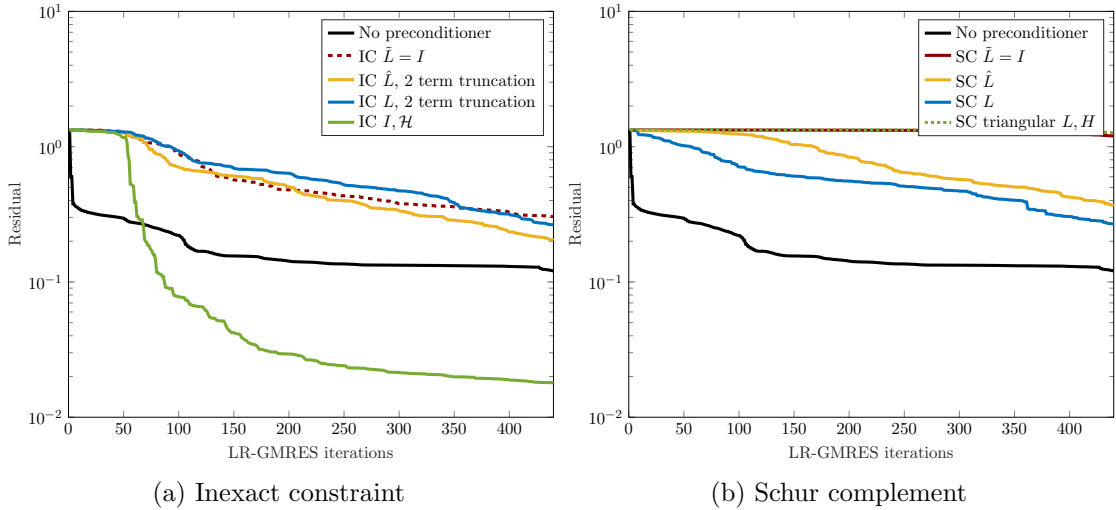

		\begin{subfigure}{.48\textwidth}
			\centering
			\resizebox{\textwidth}{!}{\input{tikzplots/ICprec440.tex}}
			\caption{Inexact constraint}
			\label{fig:IC}
		\end{subfigure}
		\begin{subfigure}{.48\textwidth}
			\centering
			\resizebox{\textwidth}{!}{\input{tikzplots/SCprec440.tex}}
			\caption{Schur complement}
			\label{fig:SC}
		\end{subfigure}
		\caption{Residual using different preconditioners for the $440 \times 440$ advection-diffusion example.}
		\label{fig:Prec}
	\end{figure}

	The preconditioners considered in Figure~\ref{fig:IC} are inexact constraint preconditioners \eqref{eq:IC_L0}, which we compare to using no preconditioner. We use $\tilde{L} = I, \hat{L} = I_{N+1} \kron I_n + C \kron I_n$, and also consider $\mathcal{P}_{I\mathcal{H}}$ from \eqref{eq:IC_IH} where $\tilde{L}=I$, and use the exact $\mathcal{H}$.
	
	We see that none of the preconditioners achieve a residual smaller than $10^{-2}$ even after $440$ iterations due to the additional restrictions of the low-rank solver (e.g. the truncation during the algorithm). The three inexact constraint preconditioners where we take $\tilde{\mathcal{H}} = 0$ exhibit very similar behaviour with the approximation $\hat{L}$ performing slightly better than the other two on the whole. The only preconditioner which achieved superior results to taking the identity, was $\mathcal{P}_{I\mathcal{H}}$ from \eqref{eq:IC_IH}, incorporating the true $\mathcal{H}$ and taking $L = I$. Despite this, the improvement occurs only after 70 iterations which for GMRES is not ideal since we must store all iterates. Even using the low-rank representation here, this becomes problematic.
	
	For Figure~\ref{fig:SC}, we experimented with a selection of Schur complement preconditioners, all of which approximate the Schur complement using the approximation \eqref{eq:SchurApprox}. For the block triangular preconditioner, we use the exact $L$ and $\mathcal{H}$ in the inverted matrix in addition to \eqref{eq:SchurApprox}.
	
	Unlike the inexact constraint preconditioners, none of the Schur complement preconditioners we consider here showed better results than using no preconditioner. Comparison with the inexact constraint preconditioners shows the block diagonal Schur complement preconditioners using $\hat{L}$ and $L$ to be comparable. Despite the block triangular preconditioner containing the true $\mathcal{H}$ it results in an ineffective choice, performing worse than all others considered.
	
\begin{figure}[!ht]
		\centering
		\resizebox{.48\textwidth}{!}{\input{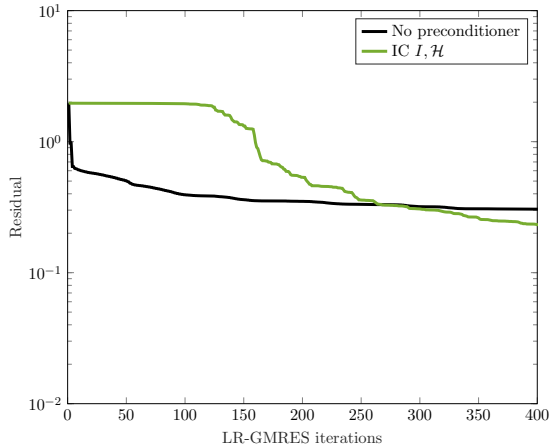}}
		\caption{Residual using the inexact constraint preconditioner for the $880 \times 880$  advection-diffusion example.}
		\label{fig:PrecIH880}
	\end{figure}

	To illustrate a larger problem size than those above, we conduct a further test using $n=20$ with the remaining setup unchanged from above. Thus the saddle point matrix is now of size $880$. In Figure~\ref{fig:PrecIH880} we compare the best performing of the above preconditioners, the inexact constraint preconditioner $\mathcal{P}_{I\mathcal{H}}$ from \eqref{eq:IC_IH} using $\tilde{L} = I$ and $\tilde{\mathcal{H}} = \mathcal{H}$.
	We see that as before, the inexact constraint preconditioner eventually results in a lower residual, though here it takes over $250$ iterations, nearly four times as many as in the $440$ system which was merely half the size. As mentioned above this is infeasible for this implementation of LR-GMRES, and hence we used no preconditioner in the numerical examples presented in Sections~\ref{sec:1DAD} and \ref{sec:Lorenz}.
	
	A possible explanation for why preconditioning is not effective here is the following. During LR-GMRES, the truncation process selects only the most important modes, e.g. the ones belonging to larger eigenvalues, ignoring the smaller ones. Therefore, the low-rank approach acts like a regularisation, and hence in some sense like a projected preconditioner.
	
	\section{Time-dependent systems}
	\label{sec:TimeDependent}
	
	Next we consider an extension of the Kronecker formulation \eqref{eq:KroneckerWeakSaddle} to the time-dependent case, allowing for time-dependent model, and observation operators, and the respective covariance matrices. The remaining assumption we must make is that the number of observations in the $i$-th timestep, $p_i$ is constant, i.e. $p_i = p$ for each $i$. With these assumptions, the linear system in \eqref{eq:KroneckerWeakSaddle} becomes
	\begin{equation}
	\label{eq:SumKroneckerWeakSaddle}
	\begin{bmatrix}
	F_1 \kron B + \sum\limits_{i=1}^{N} F_{i+1} \kron Q_{i} & 0 & I \kron I_x + \sum\limits_{i=1}^{N} C_{i} \kron M_i \\
	0 & \sum\limits_{i=0}^{N} F_{i+1} \kron R_i & \sum\limits_{i=0}^{N} F_{i+1} \kron H_i \\
	I \kron I_x + \sum\limits_{i=1}^{N} C_{i}^T \kron M_i^T & \sum\limits_{i=0}^{N} F_{i+1} \kron H_i^T & 0
	\end{bmatrix} \begin{bmatrix}
	\lambda \\ \mu \\ \delta x
	\end{bmatrix} = \begin{bmatrix}
	b \\ d \\ 0
	\end{bmatrix},
	\end{equation}
	where $F_i$ denotes the matrix with $1$ on the $i$th entry of the diagonal, and zeros elsewhere, and $C_i$ is the matrix with $-1$ on the $i$th column of the subdiagonal, and zeros elsewhere. Here $M_i$ and $H_i$ are linearisations of the model and observation operators $\mathcal{M}_i$ and $\mathcal{H}_i$ respectively about $x_i$.
	
	As in Section~\ref{sec:Kronecker}, we may use \eqref{eq:SylKronRelation} to rewrite this as the (now more general) matrix equations
	\begin{align}
	\label{eq:SumMatrixWeakSaddle}
	\begin{aligned}
	B \Lambda F_1 + \sum\limits_{i=1}^{N} Q_{i} \Lambda F_{i+1} + X + \sum\limits_{i=1}^{N} M_{i} X C_{i}^T &= \mathbbm{b} \\
	\sum\limits_{i=0}^{N} R_{i} U F_{i+1} + \sum\limits_{i=0}^{N} H_{i} X F_{i+1} &= \mathbbm{d} \\
	\Lambda + \sum\limits_{i=1}^{N} M_{i}^T \Lambda C_{i} + \sum\limits_{i=0}^{N} H_{i}^T U F_{i+1} &= 0.
	\end{aligned}
	\end{align}
	Where as before $\lambda, \delta x, b, \mu$ and $d$ are vectorised forms of the matrices $\Lambda,X,\mathbbm{b}  \in \mathbb{R}^{n \times N+1}$ and $U, \mathbbm{d}\in \mathbb{R}^{p \times N+1}$ respectively. These matrix equations must again be solved for $\Lambda$, $U$ and $X$, where $X$ is the matrix of interest.
	
	Algorithm~\ref{alg:AmultTime} is an implementation of \texttt{Amult} for the time-dependent case, explicitly writing the concatenation defined by \eqref{eq:SumMatrixWeakSaddle} in the form required for LR-GMRES. This requires linearisations of the model and observation operators at all timesteps in order to be applied.
	
	\begin{algorithm}[!ht]
		\caption{Matrix multiplication (time-dependent) (\texttt{Amult})}
		\label{alg:AmultTime}
		\begin{algorithmic}
			\REQUIRE $W_{11}, W_{12}, W_{21}, W_{22}, W_{31}, W_{32}$ 
			\ENSURE $Z_{11}, Z_{12}, Z_{21}, Z_{22}, Z_{31}, Z_{32}$ 
			\STATE $Z_{11} = [BW_{11},\quad Q_1W_{11}, \quad\ldots,\quad Q_NW_{11}, \quad  W_{31},\quad M_1W_{31}, \quad\ldots,\quad M_NW_{31}]$, 
			\STATE $Z_{12} = [F_1W_{12},\quad F_2W_{12}, \quad\ldots,\quad F_{N+1}W_{12}, \quad W_{32},\quad C_1W_{32}, \quad\ldots,\quad C_NW_{32}]$, 
			\STATE $Z_{21} = [R_0W_{21}, \quad\ldots,\quad R_{N}W_{21}, \quad H_0W_{31}, \quad\ldots,\quad H_{N}W_{31}]$,  
			\STATE $Z_{21} = [F_1W_{22}, \quad\ldots,\quad F_{N+1}W_{22}, \quad F_1W_{32}, \quad\ldots,\quad F_{N+1}W_{32}]$, 
			\STATE $Z_{31} = [W_{11}, \quad M_1^TW_{11}, \quad\ldots,\quad M_N^TW_{11}, \quad H_0^TW_{21}, \quad\ldots,\quad H_N^TW_{21}]$,  
			\STATE $Z_{32} = [W_{12},\quad C_1^TW_{12}, \quad\ldots,\quad C_N^TW_{12}, \quad F_1W_{22}, \quad\ldots,\quad F_{N+1}W_{22}]$ 
		\end{algorithmic}
	\end{algorithm}
	
	As an example, we consider the Lorenz-95 system \cite{Lorenz1996} which is both non-linear, and also chaotic rather than smoothing such as the previous example (Section~\ref{sec:1DAD}), so as to better represent real world data assimilation problems such as weather forecasting.
	
	\subsection{Lorenz-95 system}
	\label{sec:Lorenz}
	We consider the Lorenz-95 system \cite{Lorenz1996}, this is a generalisation of the three-dimensional Lorenz system \cite{Lorenz1963} to $n$ dimensions. The model is defined by a system of $n$ non-linear ordinary differential equations
	\begin{equation}
	\label{eq:Lorenz95}
	\frac{\textrm{d}x^i}{\textrm{d}t} =  -x^{i-2}x^{i-1} + x^{i-1}x^{i+1} - x^i + f,
	\end{equation}
	where $x = [x^1, x^2, \ldots, x^n]^T$ is the state of the system, and $f$ is a forcing term. It is known that for $f=8$, the Lorenz system exhibits chaotic behaviour \cite{Freitag2013,Lorenz1996}. Also noted is that for reasonably large values of $n$ (here we take $n = 40$), this choice of $f$ leads to a model which is comparable to weather forecasting models.
	
	We solve \eqref{eq:Lorenz95} using a 4th order Runge-Kutta method in order to obtain
	\begin{equation}
	x_{k+1} = \mathcal{M}_{k}(x_{k}), \quad \text{where } x_{k} = [x_k^1, x_k^2, \ldots, x_k^n]^T,
	\end{equation}
	where $\mathcal{M}_k$ is the non-linear model operator which evolves the state $x_k$ to $x_{k+1}$. As before $\mathcal{H}_k$ denotes the potentially non-linear observation operator for the state $x_k$.
	To formulate the data assimilation problem as a saddle point problem, we generate the tangent linear model, and observation operators $M_k$ and $H_k$ by linearising $\mathcal{M}_k$ and $\mathcal{H}_k$ about $x_k$.
	
	As in Section~\ref{sec:1DAD}, we compare the low-rank approximation computed using LR-GMRES, to the full-rank solution of the saddle point formulation \eqref{eq:WeakSaddle}, and the background estimate (e.g. no assimilation).
	We perform the data assimilation using an assimilation window of 200 timesteps, followed by a forecast of 1300 timesteps, where the timesteps are of size $\Delta t = 5\cdot10^{-3}$. The full-rank update is therefore $\delta x \in \mathbb{R}^{8,000}$, whilst in contrast the low-rank update $WV^T$, is such that $W \in \mathbb{R}^{40 \times r}, V \in \mathbb{R}^{200 \times r}$. Here we consider $r=20$ once more, which here requires $60\%$ of the storage, still demonstrating a significant reduction.
	
\begin{figure}[!ht]
		\begin{subfigure}{.48\textwidth}
			\centering
			\resizebox{\textwidth}{!}{
%
%
\definecolor{mycolor1}{rgb}{0.92900,0.69400,0.12500}%
\definecolor{mycolor2}{rgb}{0.49400,0.18400,0.55600}%
\definecolor{mycolor3}{rgb}{0.46600,0.67400,0.18800}%
\begin{tikzpicture}

\begin{axis}[%
width=4.359in,
height=3.651in,
at={(0.731in,0.493in)},
scale only axis,
xmin=0,
xmax=40,
xlabel style={font=\color{white!15!black}},
xlabel={x},
ymode=log,
ymin=0.0001,
ymax=10,
yminorticks=true,
ylabel style={font=\color{white!15!black}},
ylabel={Error},
axis background/.style={fill=white},
legend style={at={(0.97,0.03)}, anchor=south east, legend cell align=left, align=left, draw=white!15!black}
]
\addplot [color=mycolor1, dashed, line width=2.0pt]
  table[row sep=crcr]{%
1	0.0685889926641343\\
2	0.475737071835058\\
3	5.16274223814891\\
4	4.08193956555305\\
5	0.704729498595431\\
6	4.34447132506502\\
7	1.22647581332786\\
8	9.55189688244878\\
9	3.93635040851802\\
10	6.60886701747722\\
11	4.38209121676772\\
12	4.69784864837582\\
13	4.96201777968168\\
14	2.4747765758592\\
15	1.36153679122716\\
16	1.06835259047136\\
17	0.960539669111081\\
18	1.05581245388783\\
19	0.631446113601299\\
20	0.499556955843161\\
21	0.281834865108461\\
22	0.413672551004491\\
23	0.487691485478814\\
24	0.315098352986045\\
25	0.0993945413491758\\
26	2.23801979032544\\
27	1.58687598679372\\
28	0.473166338775627\\
29	2.67933571352564\\
30	0.116617537497394\\
31	1.88720569747536\\
32	0.733530744985326\\
33	2.10332674842699\\
34	6.69807746021849\\
35	4.95299400166489\\
36	1.63193908288492\\
37	7.47281840682778\\
38	0.199433564510612\\
39	2.51558434551802\\
40	8.7311538927116\\
};
\addlegendentry{No assimilation}

\addplot [color=mycolor2, line width=2.0pt]
  table[row sep=crcr]{%
1	0.114516535697979\\
2	0.0550959395221145\\
3	0.020648539997869\\
4	0.0256640661900995\\
5	0.233552409862613\\
6	0.00973825883857724\\
7	0.246973697439804\\
8	0.0234257840650098\\
9	0.441117292624603\\
10	0.117018134991887\\
11	0.293498827444052\\
12	0.0279971823706477\\
13	0.134054973575632\\
14	0.13338101782647\\
15	0.0582933283162682\\
16	0.0327719722283819\\
17	0.0108803036198636\\
18	0.001119064578473\\
19	0.00984359658046422\\
20	0.00741641221138017\\
21	0.00628213872472115\\
22	0.0093229592868989\\
23	0.00122597568473371\\
24	0.00499917122997793\\
25	0.000510053814500644\\
26	0.000620339964815028\\
27	0.00671451467447025\\
28	0.0411334953387063\\
29	0.00535050029813207\\
30	0.0364582522397372\\
31	0.00652670337497297\\
32	0.0472467705064651\\
33	0.0609846176164315\\
34	0.116239841088699\\
35	0.0950923947674381\\
36	0.028900615847271\\
37	0.0267534812348361\\
38	0.0727824174997885\\
39	0.00221231137544776\\
40	0.00879969371454873\\
};
\addlegendentry{Full-rank}

\addplot [color=mycolor3, line width=2.0pt]
  table[row sep=crcr]{%
1	0.110524571892527\\
2	0.0330983482335732\\
3	0.0205789436041686\\
4	0.0322445889513984\\
5	0.25266083460373\\
6	0.0210577921724196\\
7	0.251180228141052\\
8	0.0312907608720523\\
9	0.478318585474827\\
10	0.132560020720402\\
11	0.318041786022369\\
12	0.0247308825198369\\
13	0.132941122884298\\
14	0.144911201847739\\
15	0.0593786274794877\\
16	0.0310894498849892\\
17	0.0310383914821529\\
18	0.00942788204173502\\
19	0.00172712858663893\\
20	0.00654248142314318\\
21	0.0009107123996035\\
22	0.0104837565086458\\
23	0.00764292699518608\\
24	0.00611940948439617\\
25	0.00449770487803128\\
26	0.00815598130072376\\
27	0.00683580442990413\\
28	0.0208974276233649\\
29	0.00583370854191225\\
30	0.0448886028241685\\
31	0.00178247453258829\\
32	0.0451710121275825\\
33	0.0412734794972028\\
34	0.117487010363718\\
35	0.115837598751468\\
36	0.0400899896568188\\
37	0.0314882123276987\\
38	0.0781601760720766\\
39	0.013118124504365\\
40	0.00455143523989763\\
};
\addlegendentry{Low-rank}

\end{axis}
\end{tikzpicture}%
}
			\caption{Error $|x^* - x|$}
		\end{subfigure}
		\begin{subfigure}{.48\textwidth}
			\centering
			\resizebox{\textwidth}{!}{\input{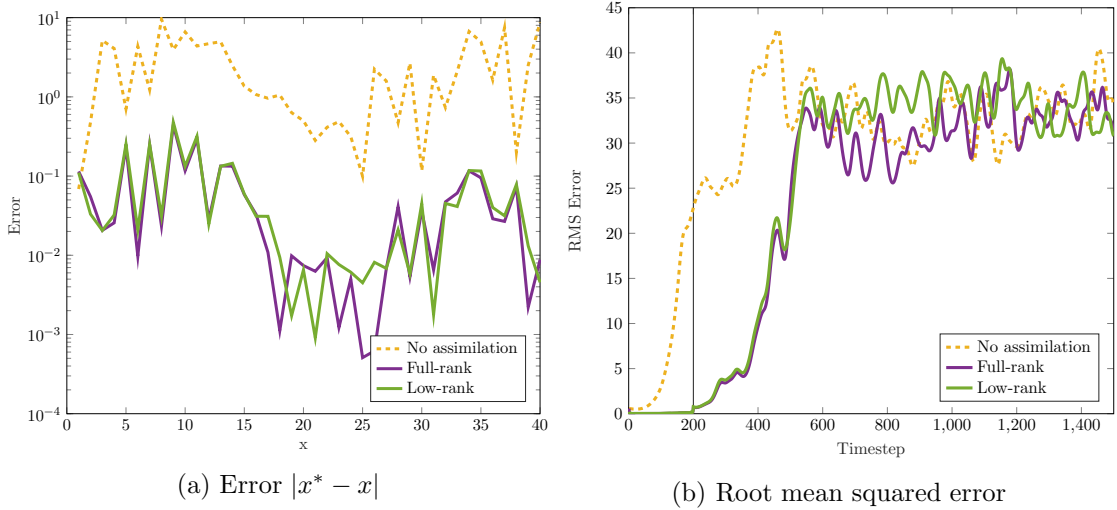}}
			\caption{Root mean squared error}
		\end{subfigure}
		\caption{Error $|x^* - x|$ for the time after the assimilation window, and root mean squared error for Lorenz-95 system with perfect observations.}
		\label{fig:L40_oN_R0_BI_E}
	\end{figure}

	\paragraph{Perfect observations}
	As with the advection-diffusion equation, let us first suppose we have perfect observations of every state in the assimilation window, we take as the background estimate $x_0^b$, a perturbed initial condition with background covariance $B = 0.1 I_{40}$, and as before, we consider a model error with covariance $Q = 10^{-4} I_{40}$.
	The error $|x^* - x|$ for the time after assimilation, and the root mean square errors for the three approaches in this example are presented in Figure~\ref{fig:L40_oN_R0_BI_E}.
	The choice of $r=20$ here results in a low-rank approximation which is very close to the full-rank solution. This is very good given that the low-rank approximation requires 40\% less storage. In the state error plot we observe small differences between solutions for the middle states, though this is still substantially smaller than the error with no assimilation. In the forecast the low-rank approximation matches the full-rank until both reach the error with no assimilation, with only small variation.
	
\begin{figure}[!ht]
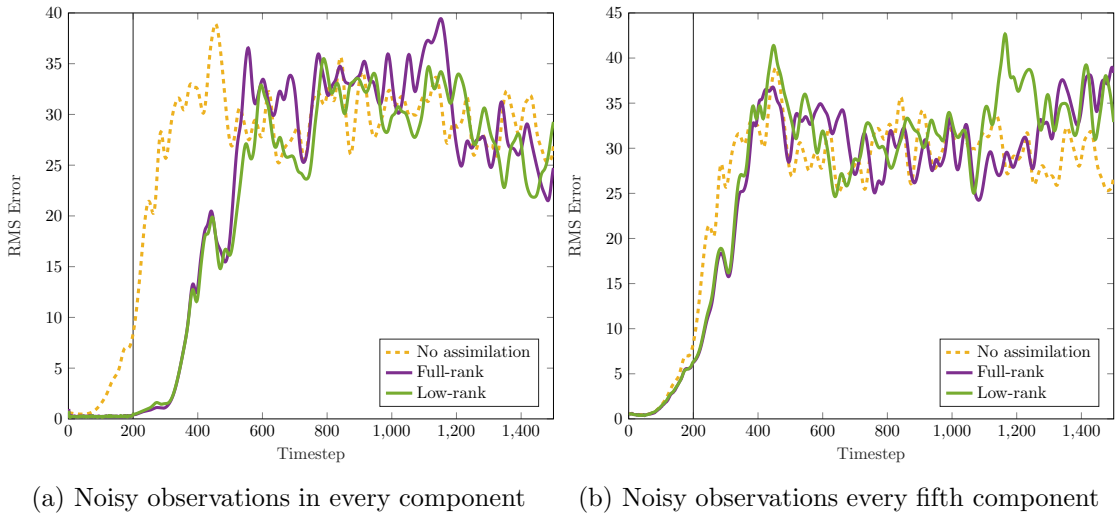

		\begin{subfigure}{.48\textwidth}
			\centering
			\resizebox{\textwidth}{!}{\input{tikzplots/L40_oN_RI_Be_RMSE.tex}}
			\caption{Noisy observations in every component}
		\end{subfigure}
		\begin{subfigure}{.48\textwidth}
			\centering
			\resizebox{\textwidth}{!}{\input{tikzplots/L40_oN5_RI_Be_RMSE.tex}}
			\caption{Noisy observations every fifth component}
		\end{subfigure}
		\caption{Root mean squared error for Lorenz-95 system with noisy, and partial observations.}
		\label{fig:L40_o_RI_Be}
	\end{figure}

	\paragraph{Noisy observations}
	We next introduce noisy observations, taking $R = 0.01I_p$ for the observation error covariance, furthermore we take as the background error covariance $B_{i,j} = 0.1\exp(\frac{-|i-j|}{50})$. In Figure~\ref{fig:L40_o_RI_Be} we consider the root mean squared errors for two different choices of observation operator: taking interpolatory observations in every component $(p = 40)$ shown on the left, and in every fifth component $(p=8)$ on the right.
	In both cases, the low-rank approximation matches the full-rank very closely until the time at which both errors are comparable to the background estimate. In this example the assimilation of noisy observations in every fifth component is similarly difficult for both approaches. 
	To achieve these very similar results using the low-rank approach, despite using just $60\%$ of the storage is very promising.
	
	\FloatBarrier
	\paragraph{150-dimensional Lorenz-95}
	Finally, we consider as a larger example, the 150 - dimensional Lorenz-95 system with an assimilation window of 150 timesteps. This gives a full-rank update $\delta x \in \mathbb{R}^{22,500}$, and we consider two different choices of low-rank, $r = 20$ requiring $27\%$ of the storage, and $r=5$ needing $7\%$. In this example we take noisy observations in each state, with covariances $B_{i,j} = 0.1\exp(\frac{-|i-j|}{50})$, $R = 0.01I_{150}$ and $Q = 10^{-4} I_{150}$.
	
\begin{figure}[!ht]
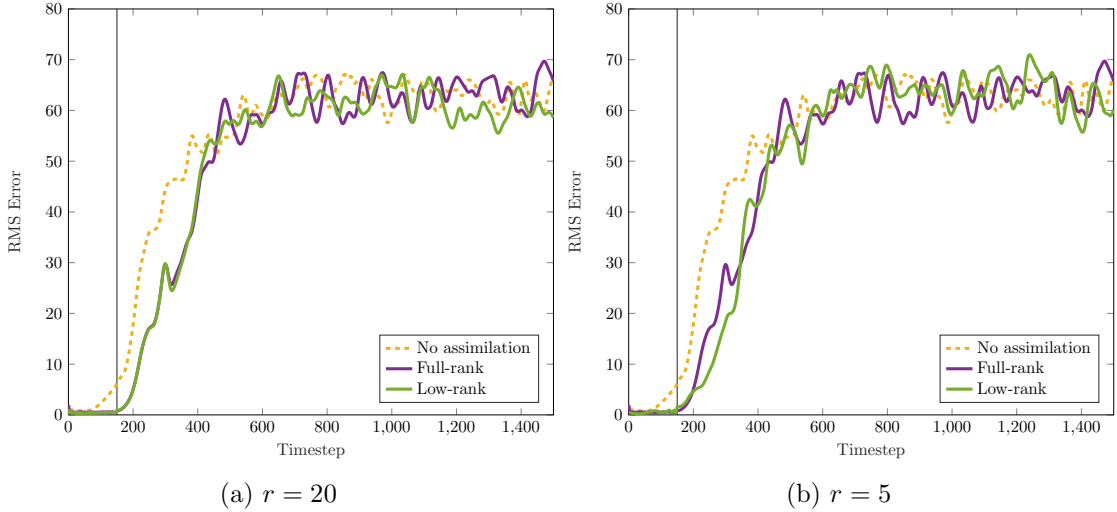

		\begin{subfigure}{.48\textwidth}
			\centering
			\resizebox{\textwidth}{!}{\input{tikzplots/LL150_oN_RI_Be_r20_RMSE.tex}}
			\caption{$r = 20$}
		\end{subfigure}
		\begin{subfigure}{.48\textwidth}
			\centering
			\resizebox{\textwidth}{!}{\input{tikzplots/LL150_oN_RI_Be_r5_RMSE.tex}}
			\caption{$r = 5$}
		\end{subfigure}
		\caption{Root mean squared error for 150-dimensional Lorenz-95 system with $r = 20$ and $r = 5$.}
		\label{fig:LL150_oN_RI_Be_RMSEr}
	\end{figure}

	These examples, shown in Figure~\ref{fig:LL150_oN_RI_Be_RMSEr} demonstrate further that a low-rank approximation performs very closely to that of the full-rank solution for small choices of $r$. Taking $r=20$ we see that as in the previous examples, the resulting approximation is nearly indistinguishable until both solutions reach the same level of error as with no assimilation.
	As before, we see the low-rank performing better for $r=5$, this is not always the case depending on the random seed as noted earlier, and is emphasised by the chaotic system sensitivity. However repeated experimentation shows that the full- and low-rank approximations are often close. Here the approximation using $r=5$ gives similar results to the full-rank approximation, despite requiring just $7\%$ of the storage.

	Table~\ref{table:Lorenzcomp} presents the storage requirements for the examples considered in this section. As with the advection-diffusion example, despite the large reduction in storage required, the experiments have shown that the low-rank approximations give similar results to the full-rank approach, which is a very good prospect.
	\begin{table}[!ht]
		\centering
		\begin{tabular}{|l|l|l|l|l|l|l|}
			\hline
			&&&& \multicolumn{2}{|l|}{\# of matrix elements in solution} & \\\cline{5-6}
			n & N & p & rank &full-rank & low-rank & storage reduction\\
			\hline
			40  & 199 & 40  & 20 & 8,000  & 4,800 &   40\% \\
			40  & 199 & 8   & 20 & 8,000  & 4,800 &   40\% \\
			150 & 149 & 150 & 20 & 22,500 & 6,000 & 73.3\% \\		
			150 & 149 & 150 & 5  & 22,500 & 1,500 & 93.3\% \\	
			\hline
		\end{tabular}
		\caption{Storage requirements for full- and low-rank methods in the Lorenz-95 examples.}
		\label{table:Lorenzcomp}
	\end{table}

	\section{Conclusions}
	
	The saddle point formulation of weak constraint four-dimensional variational data assimilation results in a large linear system which in the incremental approach is solved to determine the update $\delta x$ at every step.
	In this paper we have proposed a low-rank approach which approximates the solution to the saddle point system, with significant reductions in the storage needed. This was achieved by considering the structure of this saddle point system and using techniques from the theory of matrix equations.
	Using the existence of low-rank solutions to Sylvester equations we showed that low-rank solutions to the data assimilation problem exist under certain assumptions, with numerical experimentation demonstrating that this may be the case even when these assumptions are relaxed.
	
	We introduced a low-rank GMRES solver, considered the requirements for implementing this algorithm, and investigated several preconditioning approaches. For our examples we observed that no preconditioners were necessary, however further investigation of this may lead to new choices of preconditioners for the data assimilation setting, and new low-rank solvers for weak constraint 4D-Var.
	
	Numerical experiments have demonstrated that the low-rank approach introduced here is successful using both linear and non-linear models. In these examples we achieved close approximations to the full-rank solutions with storage requirements of up to less than $10\%$ of those needed by the full-rank approach,
	which is very promising.
	
	\FloatBarrier
	\bibliography{lowrankpaper}
	
\end{document}